\documentclass[11pt,dvipsnames]{extarticle}

\usepackage{amsfonts}
\usepackage{amsmath}
\usepackage{amssymb}
\usepackage{amsthm}
\usepackage{wasysym}
\usepackage{xcolor}
\usepackage{tikz,tikz-cd}
\usepackage{geometry}
\usepackage{physics}
\usepackage{textcomp}
\usepackage{hyperref}
\usepackage[utf8]{inputenc}
\usepackage[T1]{fontenc}
\usepackage[nameinlink, noabbrev, capitalize]{cleveref}
\usepackage{appendix}
\usepackage{mathtools}
\usepackage{xfrac}
\usepackage{nicefrac}
\usepackage{url}
\usepackage{environ}
\usepackage[hang]{footmisc}
\usepackage{algorithm2e}
\usepackage[normalem]{ulem}
\usepackage[shortlabels]{enumitem}
\usepackage{multicol}
\usepackage{xpatch}
\usepackage[font={small}]{caption}

\usepackage[backend=biber,style=alphabetic,backref=true,sorting=anyt,maxalphanames=3,minalphanames=2,maxbibnames=99,doi=false,isbn=false,url=false,eprint=false]{biblatex}

\renewbibmacro{in:}{}
\DeclareFieldFormat
[article,inbook,incollection,inproceedings,patent,thesis,unpublished]
{title}{#1}
\DefineBibliographyStrings{english}{
	backrefpage={},
	backrefpages={}
}

\DeclareFieldFormat{backrefparens}{\addperiod#1}
\xpatchbibmacro{pageref}{parens}{backrefparens}{}{}

\hypersetup{
colorlinks=true, urlcolor=NavyBlue, linkcolor=Mahogany, citecolor=ForestGreen, pdfborder={0 0 0},
}

\let\oldFootnote\footnote
\newcommand\nextToken\relax

\renewcommand\footnote[1]{%
	\oldFootnote{#1}\futurelet\nextToken\isFootnote}

\newcommand\isFootnote{%
	\ifx\footnote\nextToken\textsuperscript{,}\fi}

\mathchardef\mhyphen="2D

\newcommand{\mb}{\mathbb}

\newcommand{\mc}{\mathcal}
\newcommand{\mf}{\mathfrak}

\newcommand{\tsf}{\textsf}

\newcommand{\msf}{\mathsf}

\newcommand{\tx}{\text}
\newcommand{\ol}{\overline}
\newcommand{\ul}{\underline}

\newcommand{\spn}{\tx{\normalfont span}}

\newcommand{\wgt}{\mathrm{wt}}

\newcommand{\zcl}{\tx{\normalfont Z-cl}}
\newcommand{\zscl}{\tx{\normalfont Z*-cl}}

\newcommand{\iref}[2]{(\hyperref[#2]{\ref*{#1}.\ref*{#2}})}

\newcommand{\modulo}[1]{\,(\tx{\normalfont mod }#1)}
\newcommand{\Hilbert}{\mathrm{H}}
\newcommand{\certdeg}{\mathrm{cert\mhyphen deg}}

\newcommand{\PC}{\msf{PC}}
\newcommand{\PPC}{\msf{PPC}}

\newcommand{\email}[1]{\href{mailto:#1}{\textcolor{NavyBlue}{\texttt{#1}}}}

\newcommand{\SUt}{{\normalfont SU}\(^2\)}

\newcommand{\sqbinom}{\genfrac{[}{]}{0pt}{}}
\newcommand{\Ev}{\mathrm{Ev}}
\newcommand{\diagonal}{\mathrm{diag}}
\newcommand{\coeff}{\tx{\normalfont coeff}}
\newcommand{\Z}{\tx{\normalfont Z}}

\makeatletter
\def\thmhead@plain#1#2#3{%
	\thmname{#1}\thmnumber{\@ifnotempty{#1}{ }\@upn{#2}}%
	\thmnote{ {\the\thm@notefont\tsf{(#3)}}}}
\let\thmhead\thmhead@plain
\makeatother

\theoremstyle{plain}
\newtheorem{theorem}{Theorem}[section]
\newtheorem{fact}[theorem]{Fact}

\newtheorem{proposition}[theorem]{Proposition}

\newtheorem{lemma}[theorem]{Lemma}
\newtheorem{corollary}[theorem]{Corollary}

\newtheorem{remark}[theorem]{Remark}
\newtheorem{question}[theorem]{Question}
\newtheorem{problem}[theorem]{Problem}

\makeatletter
\NewEnviron{restatethm}[2]{%
	\protected@write\@auxout{}{%
		\string\@restatetheorem{#1}{\detokenize\expandafter{\BODY}}%
	}%
	\begin{theorem}[#2]\label{#1}\BODY\end{theorem}%
}
\newcommand{\@restatetheorem}[2]{%
	\expandafter\gdef\csname restatethm@#1\endcsname{#2}%
}
\newcommand{\restatethmnow}[2]{%
	\begingroup
	\renewcommand{\thetheorem}{\ref{#1}}%
	\begin{theorem}[#2]\csname restatethm@#1\endcsname\end{theorem}%
	\endgroup
}
\makeatother

\makeatletter
\NewEnviron{restatecor}[2]{%
	\protected@write\@auxout{}{%
		\string\@restatecorollary{#1}{\detokenize\expandafter{\BODY}}%
	}%
	\begin{corollary}[#2]\label{#1}\BODY\end{corollary}%
}
\newcommand{\@restatecorollary}[2]{%
	\expandafter\gdef\csname restatecor@#1\endcsname{#2}%
}
\newcommand{\restatecornow}[2]{%
	\begingroup
	\renewcommand{\thetheorem}{\ref{#1}}%
	\begin{corollary}[#2]\csname restatecor@#1\endcsname\end{corollary}%
	\endgroup
}
\makeatother

\makeatletter
\NewEnviron{restatelem}[2]{%
	\protected@write\@auxout{}{%
		\string\@restatelemma{#1}{\detokenize\expandafter{\BODY}}%
	}%
	\begin{lemma}[#2]\label{#1}\BODY\end{lemma}%
}
\newcommand{\@restatelemma}[2]{%
	\expandafter\gdef\csname restatelem@#1\endcsname{#2}%
}
\newcommand{\restatelemnow}[2]{%
	\begingroup
	\renewcommand{\thetheorem}{\ref{#1}}%
	\begin{lemma}[#2]\csname restatelem@#1\endcsname\end{lemma}%
	\endgroup
}
\makeatother

\title{Covering Weight-Determined Sets of Uniform Grids by Polynomials}
\author{S. Venkitesh\footnote{Institute for the Theory of Computing, Ben Gurion University of the Negev, Beersheva, Israel. Email:~\email{venkitesh.mail@gmail.com}. Webpage:~\url{https://sites.google.com/view/venkitesh}.}
}
\date{}
\begin{document}

	\maketitle

	\begin{abstract}
		We determine the minimum degree of a nonzero polynomial covering "weight-determined" subsets of a uniform grid over a field of characteristic zero.  Since such subsets are completely determined by the grid dimensions and the set of weights taken by the points, we see that the minimum degree is also determined by the set of weights.  In fact, the minimum degree can be algorithmically computed in time linear in the grid dimensions.  Previously, such a result and algorithm was known only for uniform grids that are "strictly unimodal" (Venkitesh, E-JC 2022), which implies the grid cannot be "lopsided".  Our result holds in full generality for all uniform grids.

		Our central tools are the classical algebraic objects called the "finite-degree Zariski closure" and the "affine Hilbert function".  Our technical contribution is a combinatorial characterization of both these objects for all weight-determined subsets of a uniform grid.  The minimal degree of our covering problem is precisely the threshold degree at which the behaviour of the Zariski closure changes sharply.
	\end{abstract}

%	\newpage
	\renewcommand{\baselinestretch}{0.7}\footnotesize
	\tableofcontents
	\renewcommand{\baselinestretch}{1.0}\normalsize

	\newpage

	\section{Introduction and overview}\label{sec:intro}

	\paragraph*{Notations.}  \(\mb{R}\) denotes the set of all real numbers, \(\mb{Z}\) denotes the set of all integers, \(\mb{N}\) denotes the set of all nonnegative integers, and \(\mb{Z}^+\) denotes the set of all positive integers.  For integers \(a\le b\), the interval of integers between \(a\) and \(b\) is denoted by \([a,b]\), and \([1,n]\) is also denoted by \([n]\).

	We work over the field \(\mb{R}\).  Fix \(n\in\mb{Z}^+\), and let
	\[
	G=[0,k_1-1]\times\cdots\times[0,k_n-1],
	\, k_i\ge2\text{ for all }i\in[n],
	\]
	be a uniform grid, and put \(N=\sum_{i\in[n]}(k_i-1)\).  For \(x=(x_1,\ldots,x_n)\in G\), its \tsf{weight} is \(\wgt(x)=\sum_{i\in[n]}x_i\).  A set \(A\subseteq G\) is \tsf{weight-determined} if for every \(x\in A,\,y\in G\) with \(\wgt(y)=\wgt(x)\), we have \(y\in A\). For \(E\subseteq[0,N]\), let \(\ul E\) be the set of all points of weights in \(E\); for \(E=\{j\}\), we write \(\ul j\) and call it a \tsf{layer}.  Thus the weight-determined sets are precisely the unions of layers.  We identify \(E\) with \(\ul E\) whenever convenient, and write \(\sqbinom{G}{E}=|\ul E|\) and \(\sqbinom{G}{j}=|\ul j|\).

	A classical result of de Bruijn, Tengbergen, and Kruyswijk~\cite{debruijn1951set} gives a symmetric chain decomposition of \(G\).  In particular, the layer sizes are symmetric and unimodal:
	\[
	\sqbinom{G}{0}\le\cdots\le\sqbinom{G}{\lfloor N/2\rfloor}=\sqbinom{G}{\lceil N/2\rceil}\ge\cdots\ge\sqbinom{G}{N}.
	\]
	More strongly, \(G\) is a Peck poset~\cite{proctor1980product,proctor1982representations}.

	For comparison, consider the Boolean cube \(\{0,1\}^n\).  Its weight is the Hamming weight, \(|x|\coloneqq\wgt(x)=\sum_{i\in[n]}x_i\). A subset \(A\subseteq\{0,1\}^n\) is symmetric if it is invariant under every permutation of the coordinates, that is, \((x_1,\ldots,x_n)\in A,\,\sigma\in\mf{S}_n\quad\implies\quad(x_{\sigma(1)},\ldots,x_{\sigma(n)})\in A\). The symmetric subsets of the Boolean cube are exactly its weight-determined subsets.  This equivalence does not hold for general uniform grids.
	\begin{remark}
		In general, over uniform grids, weight-determined sets and symmetric sets (sets invariant, as subsets of \(\mb{R}^n\), under permutations of coordinates) are different.  For instance, take a uniform grid \(G=[0,k_1-1]\times\cdots\times[0,k_n-1]\), where \(n\ge2\), \(k_i\ge3\) for all \(i\in[n]\), and \(k_i\ne k_j\) for some \(i,j\in[n],\,i\ne j\).  Then trivially, \(G\) is a weight-determined set that is not symmetric.  Further, the Boolean cube \(\{0,1\}^n\subseteq G\) is a symmetric set that is not weight-determined.
	\end{remark}

	\subsection{Our polynomial covering problem}

	Let \(E\subsetneq[0,N]\).  A \tsf{nontrivial polynomial cover} of \(\ul E\) is a polynomial that vanishes on \(\ul E\) but does not vanish on all of \(G\).  A \tsf{proper polynomial cover} vanishes on \(\ul E\) without vanishing identically on any layer outside \(E\).  We denote their least possible degrees by \(\PC_G(E)\) and \(\PPC_G(E)\), respectively.  The exact definitions, together with the elementary polynomial that covers precisely the layers in \(E\), are given in~\cref{sec:prelims}; see in particular Problem~\ref{prob:polynomial-covering}.

	Our objective is to determine \(\PC_G(E)\) and \(\PPC_G(E)\) for every proper weight set \(E\), and to understand the associated certifying degrees and algorithms.  The answers depend on the maximum plateau of the layer-size sequence.  Put
	\[
	p=\min\left\{N-\max_{i\in[n]}(k_i-1),\left\lfloor\frac N2\right\rfloor\right\},
	\, q=N-p,
	\]
	and, for \(i\in[0,\lfloor N/2\rfloor]\), put \(T_{N,i}=[0,i-1]\cup[N-i+1,N]\). The integer \(p\) is the first weight at which the maximum layer size is attained.  The two covering problems have the following exact answers:
	\[
	\PC_G(E)=
	\min_{\substack{0\le i\le p\\E\cup T_{N,i}\ne[0,N]}}
	\left(i+|E\setminus T_{N,i}|\right),
	\,
	\PPC_G(E)=
	|E|-\max\{i\in[0,p]:T_{N,i}\subseteq E\}.
	\]
	The first formula determines the least degree of a nontrivial cover, and the second determines the least degree of a proper cover; see~\cref{thm:PC-formula,thm:PPC-formula}.  For every nonempty proper \(E\), the certifying degree of \(\ul E\) is \(\PC_G(E)\).  The empty set is the only exception: its certifying degree is one, while both covering degrees are zero.  Moreover, the degree-\(d\) \(\Z^*\)-closure and all these parameters can be computed directly in time \(O(n+N)\); see~\cref{thm:certdeg,thm:covering-algorithms}.

	\subsection{Finite-degree \texorpdfstring{\(\Z\)- and \(\Z^*\)-closures}{Z- and Z*-closures}}

	The tool that connects the two covering problems is finite-degree closure.  For \(A\subseteq G\), the degree-\(d\) \(\Z\)-closure \(\zcl_{G,d}(A)\) is the common zero set in \(G\) of all polynomials of degree at most \(d\) that vanish on \(A\).  This is a point set in \(G\).  For a weight set \(E\subseteq[0,N]\), the degree-\(d\) \(\Z^*\)-closure is the weight set \(\zscl_{G,d}(E)=\{j\in[0,N]:\ul{j}\subseteq\zcl_{G,d}(\ul{E})\}\). Thus \(\ul{\zscl_{G,d}(E)}\) is the largest weight-determined subset of \(\zcl_{G,d}(\ul E)\).  These two closure objects have different types: \(\zcl_{G,d}(\ul E)\subseteq G\), whereas \(\zscl_{G,d}(E)\subseteq[0,N]\).

	It follows easily by definitions, that the covering parameters are precisely the two natural threshold degrees of \(\Z^*\)-closure:
	\begin{align*}
	\PC_G(E)&=\min\{d\in[0,N]:\zscl_{G,d}(E)\ne[0,N]\},\\
	\PPC_G(E)&=\min\{d\in[0,N]:\zscl_{G,d}(E)=E\}.
	\end{align*}
	This reduction is proved in Proposition~\ref{prop:covering-bridges}.  Consequently, the polynomial covering problem is reduced to determining \(\zscl_{G,d}(E)\) for every degree and every weight set.

	\subsection{Characterizing \texorpdfstring{\(\Z^*\)-closure}{Z*-closure}}

	The earlier closure and covering theory in~\cite{venkitesh-2022-covering-symmetric-sets} applies to strictly unimodal uniform grids, for which
	\[
	\sqbinom{G}{0}<\cdots<\sqbinom{G}{\lfloor N/2\rfloor}=\sqbinom{G}{\lceil N/2\rceil}>\cdots>\sqbinom{G}{N}.
	\]
	Strict unimodality can be recognized directly from the side lengths.
	\begin{theorem}[{\normalfont\cite{dhand-2014-grid-unimodal}}]
		A uniform grid \(G=[0,k_1-1]\times\cdots\times[0,k_n-1]\) is strictly unimodal if and only if
		\[
		\max_{i\in[n]}(k_i-1)\le\left\lceil\frac{1}{2}\sum_{i\in[n]}(k_i-1)\right\rceil.
		\]
	\end{theorem}
	More generally, with \(p,q\) as defined above, we determine the complete layer-size sequence:
	\[
	\sqbinom{G}{0}<\cdots<\sqbinom{G}{p}
	=\cdots=
	\sqbinom{G}{q}>\cdots>\sqbinom{G}{N}.
	\]
	Thus \([p,q]\) is the unique maximum plateau, and strict unimodality is exactly the condition \(q-p\le1\); see Theorem~\ref{thm:layer-sequence} and Corollary~\ref{cor:SU2-criterion}.

	For \(t\in[0,N]\), let \(t^\oplus\) and \(t^\ominus\) be the right and left endpoints of the maximal interval of equal layer sizes containing \(t\):
	\begin{align*}
	t^\oplus&=\max\bigg\{s\in[t,N]:\sqbinom{G}{r}=\sqbinom{G}{t}\tx{ for all }r\in[t,s]\bigg\},\\
	t^\ominus&=\min\bigg\{s\in[0,t]:\sqbinom{G}{r}=\sqbinom{G}{t}\tx{ for all }r\in[s,t]\bigg\}.
	\end{align*}
	For \(d\in[0,N]\) and \(E=\{t_1<\cdots<t_s\}\subseteq[0,N]\), define
	\[
	L_{N,d}(E)=\begin{cases}
		E&\tx{if }s\le d,\\
		\big[0,t_{s-d}^\oplus\big]\cup E\cup\big[t_{d+1}^\ominus,N\big]&\tx{if }s\ge d+1.
	\end{cases}
	\]
	Let \(L_{N,d}^{k+1}=L_{N,d}\circ L_{N,d}^k\), and set \(\ol L_{N,d}(E)=\bigcup_{k\ge1}L_{N,d}^k(E)\).  The finite-degree \(\Z^*\)-closure is exactly this stabilized iterate.
	\begin{theorem}\label{thm:Zscl-Lbar}
		Let \(G\) be a uniform grid.  For any \(d\in[0,N]\) and \(E\subseteq[0,N]\),
		\[
		\zscl_{G,d}(E)=\ol{L}_{N,d}(E).
		\]
	\end{theorem}
	This stabilized theorem has a direct, noniterative form.  Examine \(i=0,1,\ldots,\min\{d,p\}\), and let \(i_*\) be the least index satisfying
	\[
	E\cup T_{N,i_*}\ne[0,N]
	\,\text{and}\,
	|E\setminus T_{N,i_*}|\le d-i_*.
	\]
	If such an index exists, then \(\zscl_{G,d}(E)=E\cup T_{N,i_*}\); and, if no such index exists, then \(\zscl_{G,d}(E)=[0,N]\).  This is~\cref{thm:closure-direct}.  It is the form that yields directly the covering formulas and the linear-time algorithms stated above.

	\subsection{Affine Hilbert functions as the method}

	We obtain the closure theorem through affine Hilbert functions.  For a finite set \(A\subseteq\mb{R}^n\), let \(V_d(A)\) be the vector space of functions on \(A\) which admit a polynomial representation of degree at most \(d\).  Its affine Hilbert function is \(\Hilbert_d(A)=\dim V_d(A)\).  The basic criterion proved in Proposition~\ref{prop:Hilbert-closure-criterion} is
	\[
	w\in\zscl_{G,d}(E)
	\quad\Longleftrightarrow\quad
	\Hilbert_d(E\cup\{w\})=\Hilbert_d(E).
	\]
	In particular, the Hilbert function is constant on a weight set and its \(\Z^*\)-closure.
	\begin{fact}\label{fact:H-zscl}
		Let \(G\) be a uniform grid.  For any \(d\in[0,N]\) and \(E\subseteq[0,N]\), we have \(\Hilbert_d(E)=\Hilbert_d(\zscl_{G,d}(E))\).
	\end{fact}
	Therefore a complete description of the affine Hilbert functions of weight-determined sets gives a complete description of their finite-degree \(\Z^*\)-closures, and hence of the two covering degrees.  This leads naturally to the following broader question.
	\begin{question}
		Fix a field \(\mb{F}\).  Let \(G\) be a uniform grid (if \(\mb{F}\) has zero characteristic), or the Boolean cube \(\{0,1\}^n\) (if \(\mb{F}\) has positive characteristic).  For every \(E\subseteq[0,N]\) (\(N=n\) in the Boolean cube setting), characterize (combinatorially) the affine Hilbert function \(\Hilbert_d(E),\,d\in[0,N]\).
	\end{question}

	For the Boolean cube in characteristic zero, Bernasconi and Egidi introduced the following enumeration.  Let \(N\in\mb{Z}^+\), \(d\in[0,N]\), and \(E\subseteq[0,N]\), and write
	\[
	[0,d]\setminus E=\{t_\ell<\cdots<t_1\}\quad\tx{and}\quad E\setminus[0,d]=\{w_1<\cdots<w_r\}.
	\]
	The pair of sequences is called the \tsf{\((N,d)\)-BE enumeration} of \(E\).  Their theorem determines the affine Hilbert function of every symmetric set of the Boolean cube.
	\begin{theorem}[\cite{bernasconi-egidi-hilbert}]\label{thm:BE-hilbert}
		Consider the Boolean cube \(\{0,1\}^n\).  For any \(d\in[0,n]\) and \(E\subseteq[0,n]\), if \(\big(\{t_\ell<\cdots<t_1\},\{w_1<\cdots<w_r\}\big)\) is the \((n,d)\)-BE enumeration of \(E\),
		\[
		\Hilbert_d(E)=\sum_{w\in E\cap[0,d]}\binom{n}{w}+\sum_{j=1}^{\min\{\ell,r\}}\min\bigg\{\binom{n}{t_j},\binom{n}{w_j}\bigg\}.
		\]
	\end{theorem}
	They used this theorem to generalize Smolensky's result on approximation by low-degree polynomials~\cite{smolensky-low-degree}.

	There are also positive-characteristic results for special symmetric sets.  Fix a prime \(p\), let \(q\) be a power of \(p\), and define \(E_{n,a,s,q}=\{j\in[0,n]:j\equiv t\modulo{q},\tx{ for some }t\in[a,a+s-1]\}\).  We use the convention \(\binom{n}{u}=0\) whenever \(u\notin[0,n]\). Felszeghy, Heged\H{u}s, and R\'onyai~\cite{felszeghy-hegedus-ronyai-2009-complete-wide} proved the following theorem.
	\begin{theorem}[\cite{felszeghy-hegedus-ronyai-2009-complete-wide}]\label{thm:FHR}
		Fix the field \(\mb{F}_p\), where \(p\) is a prime, and let \(q\) be a power of \(p\).  Let \(a\in[0,n]\) and \(s\in[q-1]\) such that \((n-s-q)/2<a\le(n-s+q)/2\), and let \(m=\min\{a+s-1,n-a\}\).
		\begin{itemize}
			\item  If \(d\in[0,m]\), then
			\[
			\Hilbert_d(E_{n,a,s,q})=\sum_{i=0}^{\lfloor d/q\rfloor}\sum_{k=0}^{s-1}\binom{n}{d-iq-k}.
			\]
			\item  If \(d\in[m+1,n]\), then
			\[
			\Hilbert_d(E_{n,a,s,q})=\sum_{i=-\lfloor m/q\rfloor}^{\lfloor(n+s-m-1)/q\rfloor}\sum_{k=0}^{s-1}\binom{n}{m+iq-k}-\sum_{i=1}^{\lfloor(n+s-d-1)/q\rfloor}\sum_{k=0}^{s-1}\binom{n}{d+iq-k}.
			\]
		\end{itemize}
	\end{theorem}
	They used it to obtain an upper bound for set systems satisfying restricted intersection conditions, thereby proving and generalizing a conjecture of Babai and Frankl~\cite{babai-frankl-linear}.

	Our Hilbert-function theorem extends the characteristic-zero Boolean formula to every weight-determined subset of every uniform grid.  Recall that \(\sqbinom{G}{j}=|\ul j|\).
	\begin{restatethm}{thm:grid-Hilbert}{Affine Hilbert functions of general weight-determined sets}
		\phantom{}

		\noindent Let \(G\) be a uniform grid.  For any \(d\in[0,N]\) and \(E\subseteq[0,N]\), if \(\big(\{t_\ell<\cdots<t_1\},\{w_1<\cdots<w_r\}\big)\) is the \((N,d)\)-BE enumeration of \(E\), then
		\[
		\Hilbert_d(E)=\sum_{w\in E\cap[0,d]}\sqbinom{G}{w}+\sum_{j=1}^{\min\{\ell,r\}}\min\bigg\{\sqbinom{G}{t_j},\sqbinom{G}{w_j}\bigg\}.
		\]
	\end{restatethm}
	Thus~\cref{thm:grid-Hilbert} extends~\cref{thm:BE-hilbert} from the Boolean cube to uniform grids.

	\subsection{How the Hilbert function characterization is proved}\label{subsec:bottlenecks}

	The affine Hilbert function characterization requires two rank inputs that do not pass directly from the Boolean cube to larger products of chains.  We describe the two bottlenecks and the mechanisms that resolve them.

	Let \(\mb{X}=(X_1,\ldots,X_n)\), and for \(a\in\mb N\) put \(Y_i^{(a)}=X_i(X_i-1)\cdots(X_i-a+1)\).  For \(\alpha\in\mb N^n\), put \(\mb Y^{(\alpha)}=\prod_{i\in[n]}Y_i^{(\alpha_i)}\).  For \(D,E\subseteq[0,N]\), let \(\Ev_{D,E}\) be the corresponding evaluation matrix:
	\[
	\Ev_{D,E}(\alpha,\beta)=\mb{Y}^{(\alpha)}(\beta)=\alpha!\binom{\beta}{\alpha},\quad\tx{for all }\alpha\in\ul{D},\,\beta\in\ul{E}.\footnote{Note that we have \(\alpha!\coloneqq\prod_{i\in[n]}\alpha_i!\) and \(\binom{\beta}{\alpha}\coloneqq\prod_{i\in[n]}\binom{\beta_i}{\alpha_i}\), for any \(\alpha,\beta\in\mb{N}^n\).}
	\]
	The filtered falling-factorial basis gives \(\Hilbert_d(E)=\rank(\Ev_{[0,d],E})\).

	The first bottleneck is the rank of a single layer.  Wilson's incidence-matrix formula~\cite[Theorem 1]{wilson-1990-diagonal-incidence} gives the needed Boolean-cube rank, but its inductive proof does not extend directly to general uniform grids.  Instead, we use the up operator on the ranked poset \(G\).  If \((P,\le)\) is a ranked poset with rank function \(\rho\), then
		\[
		U(a)=\sum_{\substack{b\in P,\,a\le b\\\rho(b)=\rho(a)+1}}b,\quad\tx{for all }a\in P.
		\]
	A symmetric Jordan basis for \(V(G)\) makes the relevant up operators injective below the middle rank.  This resolves the first bottleneck and yields the single-layer Hilbert function.  For the symmetric Jordan basis, note that we have the classical existence results of Canfield~\cite{canfield1980sperner}, Proctor, Saks, and Sturtevant~\cite{proctor1980product}, and Proctor~\cite{proctor1982representations}, together with Srinivasan's constructive proof~\cite{srinivasan2011symmetric}.

	The second bottleneck is the following special rank statement, which follows in the Boolean cube case by induction on the number of coordinates, but does not go through in the case of larger grids (especially lopsided ones).
		\begin{restatelem}{lem:rank-special-case}{Special upper-layer rank lemma}
			Let \(G\) be a uniform grid.  If \(d\in[0,\lfloor N/2\rfloor]\) and \(\varnothing\ne E\subseteq[N-d+1,N]\), then
			\[
			\rank(\Ev_{d,E})=\rank(\Ev_{d,\min E}).
			\]
		\end{restatelem}
	For the complete proof of~\cref{thm:grid-Hilbert}, we require the stronger theorem.
		\begin{restatethm}{thm:upper-layer-rank}{Strong upper-layer rank theorem}
			Let \(G\) be a uniform grid.  If \(d\in[0,N]\) and \(\varnothing\ne E\subseteq[d,N]\), then
			\[
			\rank(\Ev_{d,E})=\rank(\Ev_{d,\min E}).
			\]
		\end{restatethm}
	The special lemma follows immediately from the strong theorem.  To prove the strong theorem, we establish the adjacent evaluation-transfer identity
\[
\Ev_{d,w+1}
=
\frac{1}{w+1-d}
\Ev_{d,w}\diagonal_w^{-1}U_{w,w+1}\diagonal_{w+1}.
\]
Iteration shows that every higher-layer block lies in the column space of the minimum-layer block, which gives the required rank equality.

The logical proof chain is therefore
\[
\begin{aligned}
\text{falling-factorial evaluation transfer}
&\longrightarrow \text{upper-layer ranks}
\longrightarrow \text{affine Hilbert functions}\\
&\longrightarrow \Z^*\text{-closures}
\longrightarrow \text{polynomial covers}.
\end{aligned}
\]

	\subsection{Related work}

	\paragraph*{Polynomial covering and finite-degree Zariski closure.}  From the viewpoint of the present covering objective, the closest predecessor is the author's earlier work~\cite{venkitesh-2022-covering-symmetric-sets}.  It studies hyperplane covers of symmetric subsets of the Boolean cube and polynomial covers of weight-determined subsets of strictly unimodal uniform grids, using finite-degree \(\Z^*\)-closure as the central tool.  The present paper retains this covering-to-closure reduction and removes the strict-unimodality restriction by determining the full maximum plateau, the \(\Z^*\)-closure for every uniform grid, and the resulting covering formulas.  Finite-degree Zariski closure was formally introduced by Nie and Wang~\cite{nie2015hilbert} in connection with the polynomial method in combinatorial geometry, while related ideas appear earlier in work on generalized Hamming weights and low-degree polynomial methods; see Wei~\cite{wei-1991-GHM}, Heijnen and Pellikaan~\cite{heijnen-pellikaan-1998-GHM-Reed-Muller}, Keevash and Sudakov~\cite{keevash-sudakov-2005-min-rank-inclusion}, and Ben-Eliezer, Hod, and Lovett~\cite{ben-eliezer-hod-lovett-2012-low-degree-polys}.

	\paragraph*{Affine Hilbert functions.}  Exact affine Hilbert functions of prescribed unions of layers determine their finite-degree \(\Z^*\)-closures.  Besides the results of Bernasconi-Egidi and Felszeghy-Heged\H{u}s-R\'onyai stated above, Heged\H{u}s and R\'onyai~\cite{hegedus-ronyai-2003-grobner-complete-uniform} characterized the reduced Gr\"obner basis of a single Boolean layer for every lexicographic order and extended the analysis to linear Sperner families in~\cite{hegedus-ronyai-2018-linear-sperner}.  Srinivasan and the author~\cite{venkitesh-2021-zariski-positive-char} characterized finite-degree Z-closures for a subclass of symmetric sets in positive characteristic.  Affine Hilbert functions also appear in Reed-Muller and polar coding theory; see Abbe, Shpilka, and Ye~\cite{abbe-shpilka-ye}.  For general introductions to vanishing ideals and affine Hilbert functions, see Cox, Little, and O'Shea~\cite{cox2015ideals}, and Nie and Wang~\cite{nie2015hilbert}.

	\paragraph*{The uniform-grid formula and arbitrary finite grids.}  The current work obtains the Hilbert function characterization in full generality of weight-determined sets and places it inside the covering-to-closure program.  Golovnev, Guo, Hatami, Nagargoje, and Yan~\cite{GGHNY-2024-hilbert} determined the minimum affine Hilbert function among all subsets of a prescribed cardinality in an arbitrary finite grid and obtained corresponding extremal bounds for ordinary finite-degree closure over finite fields.  Their problem is extremal over all subsets of a given size; here the objective requires the exact affine Hilbert function of a specified union of weight layers.

%	\paragraph{Scope.}  The exact formulas are stated for uniform grids and, by Proposition~\ref{prop:affine-normalization}, for their coordinatewise affine arithmetic-progression images with the transported index weight.  Arbitrary nonuniform grids can behave differently even when the dimensions and layer sizes agree; see Remark~\ref{rem:nonuniform-grid}.  All theorems transfer to fields of characteristic zero by~\cref{thm:char-zero-transfer}.  They do not extend unchanged to small positive characteristic; an explicit rank obstruction is given in~\cref{subsec:positive-char-obstruction}.  The precise limits of the paper are summarized in~\cref{subsec:scope-limits}.

	\paragraph*{Organization of the paper.}  In~\cref{sec:prelims}, we state the covering problem formally, introduce point closure and \(\Z^*\)-closure, and prove that the two covering degrees are closure thresholds.  In~\cref{sec:Hilbert-layer}, we determine the affine Hilbert function of one layer.  In~\cref{sec:upper-layer-transfer}, we prove the evaluation-transfer identity and the strong upper-layer rank theorem.  In~\cref{sec:Hilbert-general}, we determine the affine Hilbert functions of arbitrary unions of layers.  In~\cref{sec:layer-sequence}, we determine the exact layer sequence and its maximum plateau.  In~\cref{sec:Z-star-clo}, we characterize all finite-degree \(\Z^*\)-closures.  We then return in~\cref{sec:covering} to solve the two polynomial covering problems, determine certifying degrees, and give the direct linear-time algorithms.  We conclude with some additional remarks in~\cref{sec:further-remarks}.

\section{Preliminaries and the covering-to-closure reduction}\label{sec:prelims}\label{sec:framework}

In this section, we reduce the polynomial covering problem to finite-degree \(\Z^*\)-closure.  We keep the discussion over \(\mb{R}\); the extension to arbitrary fields of characteristic zero will be recorded at the end of the paper.

We first fix the notation for uniform grids and state the two covering parameters.  We then introduce polynomial functions and finite-degree closure, and prove that the covering degrees are exactly the two natural closure thresholds.  The remaining preliminaries record the affine Hilbert and rank machinery used to compute that closure.

\subsubsection*{Uniform grids, layers, and weight-determined sets}

Fix \(n\in\mb{Z}^+\) and a uniform grid \(G=[0,k_1-1]\times\cdots\times[0,k_n-1]\) for \(k_i\ge2\ \text{ for all }i\in[n]\), and write \(N=\sum_{i\in[n]}(k_i-1), \, \top=(k_1-1,\ldots,k_n-1)\). For every \(\alpha=(\alpha_1,\ldots,\alpha_n)\in\mb{N}^n\), put \(\wgt(\alpha)=\sum_{i\in[n]}\alpha_i\). For \(x\in G\), this is the weight function used in the Introduction.  For \(w\in[0,N]\), the \tsf{weight-\(w\) layer} is \(\ul{w}=\{x\in G:\wgt(x)=w\}\). More generally, for \(E\subseteq[0,N]\), put \(\ul{E}=\bigcup_{w\in E}\ul{w}\). Thus the weight-determined subsets of \(G\) are precisely the sets \(\ul{E}\).  We will freely identify \(E\) with \(\ul{E}\) whenever the type of the object is clear.  We write \(\sqbinom{G}{E}=|\ul{E}|, \, \sqbinom{G}{w}=|\ul{w}|\). The reflection \(\rho:G\longrightarrow G, \, \rho(x)=\top-x\), is a bijection from \(\ul{w}\) to \(\ul{N-w}\).  In particular, \(\sqbinom{G}{w}=\sqbinom{G}{N-w} \,\text{for all }w\in[0,N]\).

We use the coordinatewise partial order on \(\mb{N}^n\), and hence on \(G\): for \(a,b\in\mb{N}^n\), \(a\le b \quad\Longleftrightarrow\quad a_i\le b_i\ \text{ for all }i\in[n]\). For \(a\in G\), define \(\Delta(a)=\{b\in G:b\le a\}, \, \nabla(a)=\{b\in G:b\ge a\}\). For \(A\subseteq G\), put \(\Delta(A)=\bigcup_{a\in A}\Delta(a)\), and define \(\nabla(A)\) similarly.  A set \(D\subseteq G\) is a \tsf{downset} if \(\Delta(D)=D\).  The grid is a lattice: for \(a,b\in G\), \(a\wedge b=(\min\{a_i,b_i\}:i\in[n]), \, a\vee b=(\max\{a_i,b_i\}:i\in[n])\). We denote the zero vector by \(0^n\), and the \(i\)-th standard unit vector by \(e_i\).

We may now state the problem that drives the paper.  Let \(\mb{X}=(X_1,\ldots,X_n)\), and let \(\mb{R}[\mb{X}]\) be the polynomial ring in these indeterminates.

\subsubsection*{Nontrivial and proper polynomial covers}

Recall that, for a polynomial \(P(\mb X)\), we write \(\mc Z(P)=\{x\in G:P(x)=0\}\). Let \(E\subsetneq[0,N]\).  We say that \(P\in\mb R[\mb X]\) is
\begin{itemize}
\item a \tsf{nontrivial polynomial cover} of \(\ul E\) if \(\ul E\subseteq\mc Z(P)\ne G\);
\item a \tsf{proper polynomial cover} of \(\ul E\) if \(\ul E\subseteq\mc Z(P) \quad\text{and}\quad \ul j\not\subseteq\mc Z(P) \quad\text{for every }j\in[0,N]\setminus E\).
\end{itemize}
Thus a proper cover may vanish at some points of a layer outside \(E\), but it does not cover that entire layer.

Let \(\PC_G(E)\) and \(\PPC_G(E)\) be the minimum degrees of a nontrivial polynomial cover and a proper polynomial cover of \(\ul E\), respectively.  These quantities are well-defined.  Indeed,
\begin{equation}\label{eq:weight-cover-polynomial}
P_E(\mb X)=
\prod_{t\in E}
\left(\sum_{r\in[n]}X_r-t\right)
\end{equation}
vanishes exactly on the layers indexed by \(E\).  Indeed, if \(x\in\ul w\), then \(P_E(x)=\prod_{t\in E}(w-t)\), which is zero precisely when \(w\in E\).  When \(E=\emptyset\), the empty product is the constant polynomial \(1\).  Consequently, \(0\le\PC_G(E)\le\PPC_G(E)\le|E|\le N\).

\begin{problem}\label{prob:polynomial-covering}
Determine \(\PC_G(E)\) and \(\PPC_G(E)\), for every proper weight set \(E\subsetneq[0,N]\).
\end{problem}

To solve Problem~\ref{prob:polynomial-covering}, we pass from individual covering polynomials to the space of all low-degree polynomial functions.  The dimensions of these spaces lead to the closure criterion used below.

\subsubsection*{Polynomial functions and affine Hilbert functions}

For \(\alpha\in\mb{N}^n\), write \(\mb{X}^\alpha=\prod_{i\in[n]}X_i^{\alpha_i}\). For \(P(\mb{X})\in\mb{R}[\mb{X}]\), let \(\coeff(\alpha,P)\) denote the coefficient of \(\mb{X}^\alpha\) in \(P\).  For a finite set \(A\subseteq\mb{R}^n\), let \(V(A)\) be the vector space of all functions \(A\to\mb{R}\).  For \(d\in\mb{N}\), let \(V_d(A)\subseteq V(A)\) be the subspace of functions which admit a polynomial representation of total degree at most \(d\).  The \tsf{affine Hilbert function} of \(A\) is \(\Hilbert_d(A)=\dim V_d(A)\). For a weight set \(E\subseteq[0,N]\), we use the abbreviations \(V(E)=V(\ul{E}), \, V_d(E)=V_d(\ul{E}), \, \Hilbert_d(E)=\Hilbert_d(\ul{E})\). When \(E=\{w\}\), we write \(V(w)\), \(V_d(w)\), and \(\Hilbert_d(w)\).

\subsubsection*{Finite-grid interpolation and degree-filtered bases}

The closure criterion below compares dimensions of low-degree polynomial-function spaces.  We therefore need a basis that respects the total-degree filtration.  The following elementary form of finite-grid interpolation supplies such a basis.  The proof follows immediately from~\cite{alon-1999-combinatorial} and Gaussian elimination; so we omit the proof.

\begin{proposition}[Finite-grid interpolation and filtered triangular bases]\label{prop:filtered-triangular-basis}\label{cor:CN}
Let \(m_1,\ldots,m_n\in\mb{N}\), and let
\[
B=[0,m_1]\times\cdots\times[0,m_n]\subseteq\mb{N}^n,
\,
a\in\mb{R}^n,
\,
H=a+B.
\]
Then the following hold.
\begin{enumerate}[{\normalfont(a)}]
\item The restrictions of the monomials \(\{\mb{X}^\alpha:\alpha\in B\}\) forms a basis of \(V(H)\).
\item Suppose that \(\{f_\alpha(\mb{X}):\alpha\in B\}\) satisfies
\[
\coeff(\alpha,f_\alpha)\ne0,
\,
f_\alpha(\mb{X})\in
\spn\{\mb{X}^\beta:\beta\in\mb{N}^n,\ \beta\le\alpha\}
\quad\text{for every }\alpha\in B.
\]
Then the restrictions of the \(f_\alpha\)'s form a basis of \(V(H)\).
\item Under the hypotheses of Item~(b), for every \(d\in\mb{N}\),
\[
V_d(H)=
\spn\{\mb{X}^\alpha|_H:\alpha\in B,\ \wgt(\alpha)\le d\}
=
\spn\{f_\alpha|_H:\alpha\in B,\ \wgt(\alpha)\le d\}.
\]
\end{enumerate}
\end{proposition}
%
%\begin{proof}
%For each \(b=(b_1,\ldots,b_n)\in H\), the tensor-product Lagrange polynomial
%\[
%\ell_b(\mb{X})=
%\prod_{i\in[n]}
%\prod_{\substack{t\in a_i+[0,m_i]\\t\ne b_i}}
%\frac{X_i-t}{b_i-t}
%\]
%has degree at most \(m_i\) in \(X_i\), and satisfies \(\ell_b(c)=1\) if \(c=b\), and \(\ell_b(c)=0\) for every other \(c\in H\).  Hence the monomials in Item~(a) span \(V(H)\).  Their number is \(|B|=|H|=\dim V(H)\), so they form a basis.
%
%Choose any linear extension of the coordinatewise order on \(B\).  Relative to the monomial basis from Item~(a), the coefficient matrix of the family in Item~(b) is triangular and has nonzero diagonal.  It is therefore invertible, proving Item~(b).
%
%It remains to keep track of the degree.  Every polynomial of degree at most \(d\) may be reduced, successively in the variables, modulo the monic polynomials
%\[
%\prod_{t\in a_i+[0,m_i]}(X_i-t),
%\, i\in[n],
%\]
%without increasing its total degree.  Its restriction to \(H\) is therefore represented by a linear combination of the monomials in Item~(c).  Finally, if \(\beta\le\alpha\), then \(\wgt(\beta)\le\wgt(\alpha)\).  Thus the triangular change of basis in Item~(b) preserves the filtration by total degree, and Item~(c) follows.
%\end{proof}

\subsubsection*{Falling factorials and evaluation matrices}

To turn affine Hilbert dimensions into matrix ranks, we use a triangular basis whose evaluations are governed by the product order.  For \(r\in\mb{N}\), define \(Y_i^{(r)}=X_i(X_i-1)\cdots(X_i-r+1), \, Y_i^{(0)}=1\). For \(\alpha\in\mb{N}^n\), put \(\mb{Y}^{(\alpha)}= \prod_{i\in[n]}Y_i^{(\alpha_i)}\). For \(\alpha,\beta\in\mb{N}^n\), we use the multi-index notation
\[
\alpha!=\prod_{i\in[n]}\alpha_i!,
\,
\binom{\beta}{\alpha}=\prod_{i\in[n]}\binom{\beta_i}{\alpha_i}.
\]
Then \(\mb{Y}^{(\alpha)}(\beta) = \alpha!\binom{\beta}{\alpha}\). In particular, this value is nonzero exactly when \(\alpha\le\beta\).  If, for a polynomial \(P\), we write \(\mc{Z}(P)=\{x\in G:P(x)=0\}\), then \(\mc{Z}(\mb{Y}^{(\alpha)})=G\setminus\nabla(\alpha) \,\text{for every }\alpha\in G\).

Since \(\mb{Y}^{(\alpha)}\) has leading monomial \(\mb{X}^\alpha\) and only componentwise smaller monomials below it, \cref{prop:filtered-triangular-basis} gives, for every \(d\in[0,N]\), that \(\{\mb{Y}^{(\alpha)}|_G:\alpha\in G,\ \wgt(\alpha)\le d\}\) is a basis of \(V_d(G)\). For \(D,E\subseteq[0,N]\), define the evaluation matrix \(\Ev_{D,E}\in M_{\ul{D}\times\ul{E}}(\mb{R})\) by
\[
\Ev_{D,E}(\alpha,\beta)
=
\mb{Y}^{(\alpha)}(\beta)
=
\alpha!\binom{\beta}{\alpha},
\,
\alpha\in\ul{D},\ \beta\in\ul{E}.
\]
Also put \(\diagonal_D= \diagonal(\alpha!:\alpha\in\ul{D})\). The preceding basis statement gives the rank identity
\[
\Hilbert_d(E)=\rank(\Ev_{[0,d],E})
\,
\text{for every }d\in[0,N]\text{ and }E\subseteq[0,N].
\]

\subsubsection*{Point closure and weight-set \texorpdfstring{\(\Z^*\)-closure}{Z*-closure}}

The covering problem is now converted into a vanishing question: which points, and which entire weight layers, are forced by all degree-at-most-\(d\) polynomials that vanish on a given set?  For \(A\subseteq G\) and \(d\in\mb{N}\), let \(\mc{I}_{G,d}(A)= \{P(\mb{X})\in\mb{R}[\mb{X}]:\deg P\le d,\ P(a)=0\text{ for every }a\in A\}\). The \tsf{degree-\(d\) Zariski closure}, or \tsf{degree-\(d\) Z-closure}, of \(A\) is the point set
\[
\zcl_{G,d}(A)=
\{x\in G:P(x)=0\text{ for every }P\in\mc{I}_{G,d}(A)\}
\subseteq G.
\]
For a weight set \(E\subseteq[0,N]\), its \tsf{degree-\(d\) Z*-closure} is the weight set \(\zscl_{G,d}(E)= \{w\in[0,N]:\ul{w}\subseteq\zcl_{G,d}(\ul{E})\} \subseteq[0,N]\). Equivalently, \(\ul{\zscl_{G,d}(E)}\) is the largest weight-determined subset of the point set \(\zcl_{G,d}(\ul{E})\).  We will keep these two types distinct: \(\zcl_{G,d}(A)\) is a subset of \(G\), while \(\zscl_{G,d}(E)\) is a subset of \([0,N]\).

\begin{proposition}[Basic closure properties]\label{prop:zstar-basic}
For every \(d\in\mb{N}\), the maps
\[
\zcl_{G,d}:2^G\longrightarrow2^G
\,\text{and}\,
\zscl_{G,d}:2^{[0,N]}\longrightarrow2^{[0,N]}
\]
are closure operators.  Moreover, for all \(d,d'\in\mb{N}\) with \(d\le d'\), all \(A\subseteq G\), and all \(E\subseteq[0,N]\),
\[
\zcl_{G,d'}(A)\subseteq\zcl_{G,d}(A)
\,\text{and}\,
\zscl_{G,d'}(E)\subseteq\zscl_{G,d}(E).
\]
\end{proposition}

\begin{proof}
For point closure, extensivity and monotonicity follow directly from the definition.  Further, a polynomial of degree at most \(d\) vanishes on \(A\) if and only if it vanishes on \(\zcl_{G,d}(A)\).  Hence point closure is idempotent.  Increasing the degree only increases the collection of polynomials whose common zero set is taken, which proves the displayed inclusion for point closures.

Now let \(F=\zscl_{G,d}(E)\).  By definition, \(\ul{E}\subseteq\ul{F}\subseteq\zcl_{G,d}(\ul{E})\). Monotonicity of point closure gives
\[
\zcl_{G,d}(\ul{E})
\subseteq
\zcl_{G,d}(\ul{F})
\subseteq
\zcl_{G,d}(\zcl_{G,d}(\ul{E}))
=
\zcl_{G,d}(\ul{E}).
\]
Thus the two point closures are equal, and taking their largest weight-determined subsets gives \(\zscl_{G,d}(F)=F\).  Extensivity, monotonicity, and the degree inclusion for \(\zscl_{G,d}\) follow in the same way from the corresponding properties of point closure.
\end{proof}

\begin{proposition}[Hilbert function criterion]\label{prop:Hilbert-closure-criterion}
Let \(A\subseteq B\subseteq G\), and let \(d\in\mb{N}\).  Then
\[
\Hilbert_d(A)=\Hilbert_d(B)
\quad\Longleftrightarrow\quad
B\subseteq\zcl_{G,d}(A).
\]
In particular, for every \(E\subseteq[0,N]\) and \(w\in[0,N]\),
\[
w\in\zscl_{G,d}(E)
\quad\Longleftrightarrow\quad
\Hilbert_d(E\cup\{w\})=\Hilbert_d(E),
\]
and
\[
\Hilbert_d\bigl(\zscl_{G,d}(E)\bigr)=\Hilbert_d(E).
\]
\end{proposition}

\begin{proof}
For \(C\subseteq G\), let \(\operatorname{res}_{C,d}:V_d(G)\longrightarrow V(C)\) be the restriction map.  Its image is exactly \(V_d(C)\), so its rank is \(\Hilbert_d(C)\); its kernel consists of the degree-at-most-\(d\) polynomial functions on \(G\) which vanish on \(C\).  Since \(A\subseteq B\), we have \(\ker(\operatorname{res}_{B,d})\subseteq\ker(\operatorname{res}_{A,d})\).  The two restriction maps have the same finite-dimensional domain.  Thus their ranks are equal if and only if their kernels are equal.  The latter equality holds if and only if every degree-at-most-\(d\) polynomial function which vanishes on \(A\) also vanishes on \(B\), which is precisely the condition \(B\subseteq\zcl_{G,d}(A)\).

Taking \(A=\ul{E}\) and \(B=\ul{E\cup\{w\}}\) gives the first consequence.  For the last one, put \(F=\zscl_{G,d}(E)\).  Then \(\ul{E}\subseteq\ul{F}\subseteq\zcl_{G,d}(\ul{E})\), so the first part gives \(\Hilbert_d(F)=\Hilbert_d(E)\).
\end{proof}

\subsubsection*{Closure-to-covering bridge theorems}

The following proposition completes the reduction of Problem~\ref{prob:polynomial-covering} to finite-degree \(\Z^*\)-closure: the two covering parameters are its two natural threshold degrees.

\begin{proposition}[Covering and closure]\label{prop:covering-bridges}
For every \(E\subsetneq[0,N]\),
\begin{align}
\PC_G(E)
&=
\min\{d\in[0,N]:\zscl_{G,d}(E)\ne[0,N]\},
\label{eq:PC-bridge}\\
\PPC_G(E)
&=
\min\{d\in[0,N]:\zscl_{G,d}(E)=E\}.
\label{eq:PPC-bridge}
\end{align}
\end{proposition}

\begin{proof}
Fix \(d\in[0,N]\).  A nontrivial polynomial cover of degree at most \(d\) exists precisely when there is a polynomial of degree at most \(d\) which vanishes on \(\ul E\) but is not the zero function on \(G\).  This is equivalent to \(\zcl_{G,d}(\ul E)\ne G\), and hence to \(\zscl_{G,d}(E)\ne[0,N]\).  Taking the least such degree proves~\eqref{eq:PC-bridge}.

Suppose next that \(P\) is a proper cover of degree at most \(d\).  For every \(j\notin E\), choose \(a_j\in\ul j\) with \(P(a_j)\ne0\).  Since \(P\) vanishes on \(\ul E\), we have \(a_j\notin\zcl_{G,d}(\ul E)\). Thus no layer outside \(E\) is contained in the point closure, and therefore \(\zscl_{G,d}(E)=E\).

Conversely, suppose that \(\zscl_{G,d}(E)=E\).  For every \(j\notin E\), choose \(a_j\in\ul j\setminus\zcl_{G,d}(\ul E)\). By the definition of point closure, there is a polynomial \(P_j\) of degree at most \(d\) such that \(P_j|_{\ul E}=0, \, P_j(a_j)\ne0\). Let the weights outside \(E\) be \(j_1,\ldots,j_m\).  For \(\beta=(\beta_1,\ldots,\beta_m)\in\mb R^m\), put \(P_\beta=\sum_{r=1}^m\beta_rP_{j_r}\). For each \(s\in[m]\), the condition \(P_\beta(a_{j_s})=0\) defines a proper linear hyperplane in \(\mb R^m\), since the coefficient of \(\beta_s\) in this linear equation is \(P_{j_s}(a_{j_s})\ne0\).  A finite union of proper hyperplanes cannot cover \(\mb R^m\): the product of their nonzero defining linear forms is a nonzero polynomial, and a nonzero polynomial over the infinite field \(\mb R\) cannot vanish everywhere.  We may therefore choose \(\beta\) such that \(P_\beta(a_j)\ne0 \quad\text{for every }j\notin E\). The polynomial \(P_\beta\) has degree at most \(d\), vanishes on \(\ul E\), and does not cover any entire layer outside \(E\).  It is a proper cover.  Taking the least such degree proves~\eqref{eq:PPC-bridge}.
\end{proof}

Using~\cref{thm:Zscl-Lbar}, the same identities may be written as \(\PC_G(E)=\min\{d:\ol L_{N,d}(E)\ne[0,N]\}, \, \PPC_G(E)=\min\{d:\ol L_{N,d}(E)=E\}\). Thus Problem~\ref{prob:polynomial-covering} is solved once the sets \(\zscl_{G,d}(E)\) are known.  By \cref{prop:Hilbert-closure-criterion}, this requires the affine Hilbert functions of unions of layers.  The remaining preliminaries introduce the up operator tools used in the single-layer calculation.

\subsubsection*{Symmetric Jordan bases and up operators}

The first Hilbert function calculation concerns a single layer.  For this calculation, we use the up operator on the product of chains and a symmetric Jordan basis.  For \(x\in G\), let \(\mathbf{1}_x\in V(G)\) be the point-indicator function of \(x\).  We identify \(V(\ul{w})\) with the subspace of \(V(G)\) spanned by \(\{\mathbf{1}_x:x\in\ul{w}\}\), using extension by zero outside the layer.  Define the \tsf{up operator} \(U:V(G)\to V(G)\) on this basis by \(U(\mathbf{1}_x)= \sum_{\substack{y\in G\\x<y,\ \wgt(y)=\wgt(x)+1}} \mathbf{1}_y\). For \(d\in[0,N-1]\), define \(U_{d,d+1}\in M_{\ul{d}\times\ul{d+1}}(\mb{R})\) by \(U_{d,d+1}(\alpha,\beta)=1\) for \(\alpha\le\beta\), and \(0\) otherwise. Thus the matrix of \(U:V(\ul{d})\to V(\ul{d+1})\), with respect to the point-indicator bases, is \(U_{d,d+1}^{\mathsf{t}}\).

For \(0\le r\le\lfloor N/2\rfloor\), a sequence of nonzero homogeneous vectors \(v_r,v_{r+1},\ldots,v_{N-r}, \, v_j\in V(\ul{j})\), is a \tsf{symmetric Jordan chain} for \(U\) if \(Uv_j=v_{j+1} \quad(r\le j<N-r), \, Uv_{N-r}=0\). A basis of \(V(G)\) which is a disjoint union of symmetric Jordan chains is called a \tsf{symmetric Jordan basis}.  A classical theorem says that such a basis exists for every product of finite chains; we will state the exact consequence that we need in \cref{sec:Hilbert-layer}.

The reduction is now complete: covering degrees are closure thresholds, closure membership is tested by affine Hilbert functions, and the next two sections compute those ranks.  We begin with a single layer.

	\section{Affine Hilbert function of a single layer}\label{sec:Hilbert-layer}

By \cref{prop:Hilbert-closure-criterion}, deciding whether a new layer enters finite-degree \(\Z^*\)-closure amounts to deciding whether adjoining that layer changes the affine Hilbert dimension.  The first Hilbert input is therefore the Hilbert function of a single layer \(\ul{w}\).  We determine it in this section; the two cases in the statement below are both necessary.

\begin{theorem}[Affine Hilbert function of a single layer]\label{thm:single-layer-Hilbert}
    Let \(G\) be a uniform grid and let \(d,w\in[0,N]\).  Then
    \[
    \Hilbert_d(w)=
    \begin{cases}
    \sqbinom{G}{w},&d\ge w,\\[2mm]
    \min\left\{\sqbinom{G}{d},\sqbinom{G}{w}\right\},&d<w.
    \end{cases}
    \]
    Equivalently,
    \[
    \Hilbert_d(w)=
    \min\left\{
    \sqbinom{G}{\min\{d,w\}},
    \sqbinom{G}{w}
    \right\}.
    \]
\end{theorem}

This is the first rank calculation required by the closure problem.  The proof uses the falling-factorial evaluation matrices from \cref{sec:framework} and a symmetric Jordan basis for the up operator.  We first record the precise consequences of the symmetric Jordan basis that we need.

\subsection{Up operator ranks from symmetric Jordan bases}

Several proofs of the existence of a symmetric Jordan basis for a product of chains are known; see Canfield~\cite{canfield1980sperner}, Proctor, Saks, and Sturtevant~\cite{proctor1980product}, Proctor~\cite{proctor1982representations}, and Srinivasan~\cite{srinivasan2011symmetric}.  Srinivasan gives a constructive proof.  We use only the rank consequences of such a basis.

\begin{proposition}[Consequences of a symmetric Jordan basis]\label{prop:up-matrix}
    Let \(G\) be a uniform grid.
    \begin{enumerate}[{\normalfont(a)}]
        \item For every \(d\in[0,N-1]\),
        \[
        \rank(U_{d,d+1})=
        \min\left\{\sqbinom{G}{d},\sqbinom{G}{d+1}\right\}.
        \]
        Thus the up operator is injective when \(d<N/2\), and surjective when \(d\ge\lfloor N/2\rfloor\).

        \item If \(0\le d<w\le\left\lfloor\frac N2\right\rfloor\), then \(U_{d,d+1}U_{d+1,d+2}\cdots U_{w-1,w}\) has row rank \(\sqbinom{G}{d}\).

        \item The layer sizes are symmetric and unimodal, that is, \(\sqbinom{G}{j}=\sqbinom{G}{N-j}\) for every \(j\in[0,N]\), and \(\sqbinom{G}{0}\le\sqbinom{G}{1}\le\cdots\le\sqbinom{G}{\lfloor N/2\rfloor}\).
    \end{enumerate}
\end{proposition}

\begin{proof}
    Fix a symmetric Jordan basis for the up operator \(U\).  Every Jordan chain has the form \(v_r,v_{r+1},\ldots,v_{N-r}, \, Uv_j=v_{j+1}\). If \(d<N/2\), then a chain containing a vector of weight \(d\) starts at some weight \(r\le d\), and its last weight is \(N-r\ge N-d>d\).  Thus the chain also contains a vector of weight \(d+1\), and \(U:V(\ul{d})\to V(\ul{d+1})\) is injective.  If \(d\ge\lfloor N/2\rfloor\), then a chain containing a vector of weight \(d+1\) starts at some weight \(r\le N-d-1\le d\).  Hence the chain also contains its predecessor of weight \(d\), and the same map is surjective.  Since the matrix of this restriction in the point-indicator bases is \(U_{d,d+1}^{\mathsf t}\), part~(a) follows.

    Now suppose \(d<w\le\lfloor N/2\rfloor\).  A chain containing a vector of weight \(d\) starts at some weight \(r\le d\) and ends at weight \(N-r\ge N-d>w\).  It therefore continues at least through weight \(w\).  Hence the composition \(U^{w-d}:V(\ul{d})\longrightarrow V(\ul{w})\) is injective in the symmetric Jordan basis.  The transpose of its matrix is \(U_{d,d+1}U_{d+1,d+2}\cdots U_{w-1,w}\), which proves part~(b).  Notice that this uses the Jordan chains for the whole composition; full rank of the adjacent rectangular matrices alone would not imply full rank of their product.

    Finally, a symmetric Jordan chain contains a vector of weight \(j\) if and only if it contains a vector of weight \(N-j\).  This gives symmetry of the layer sizes.  The injectivity proved above gives their lower-half unimodality, proving part~(c).
\end{proof}

\subsection{Consecutive spanning and evaluation factorization}

The up operator ranks must now be transferred to the falling-factorial evaluation matrices that compute the Hilbert function.  For a single layer, we abbreviate \(\diagonal_{\{w\}}\) by \(\diagonal_w\).  Thus \(\Ev_{w,w}=\diagonal_w\), and this matrix is invertible.

The next elementary identity is the link between the up operator and the evaluation matrices: on a fixed higher layer, it rewrites a falling factorial as a linear combination of falling factorials of the next exact degree.

\begin{lemma}\label{lem:consec-span}
    Let \(d,w\in[0,N]\) with \(d<w\), and let \(\alpha\in\ul{d}\).  Then, as functions on \(\ul{w}\),
    \[
    \mb{Y}^{(\alpha)}=
    \frac{1}{w-d}
    \sum_{\substack{\beta\in\ul{d+1}\\\beta\ge\alpha}}
    \mb{Y}^{(\beta)}.
    \]
\end{lemma}

\begin{proof}
    Put \(I_\alpha=\{i\in[n]:\alpha_i<k_i-1\}\). Since \(\mb{Y}^{(\alpha+e_i)}= \mb{Y}^{(\alpha)}(X_i-\alpha_i)\), we have
    \[
    \sum_{i\in I_\alpha}\mb{Y}^{(\alpha+e_i)}
    =
    \mb{Y}^{(\alpha)}
    \sum_{i\in I_\alpha}(X_i-\alpha_i).
    \]
    Let \(x\in\ul{w}\).  If \(\mb{Y}^{(\alpha)}(x)=0\), then every summand on the left also vanishes, since it contains \(\mb{Y}^{(\alpha)}\) as a factor.  Otherwise \(x\ge\alpha\).  For every \(i\notin I_\alpha\), we then have \(x_i=\alpha_i=k_i-1\), and consequently \(\sum_{i\in I_\alpha}(x_i-\alpha_i) =\wgt(x)-\wgt(\alpha)=w-d\). Dividing by \(w-d\) proves the identity.
\end{proof}

Repeated application gives the exact-degree spanning statement needed to read the single-layer Hilbert function from one evaluation block.

\begin{corollary}\label{cor:full-span}
    Let \(d,w\in[0,N]\) with \(d<w\).  Then
    \[
    V_d(w)=
    \spn\{\mb{Y}^{(\alpha)}|_{\ul{w}}:\alpha\in\ul{d}\}.
    \]
    In particular, we have \(\Hilbert_d(w)=\rank(\Ev_{d,w})\).
\end{corollary}

\begin{proof}
    By \cref{prop:filtered-triangular-basis}, the restrictions of the falling factorials of weights at most \(d\) span \(V_d(w)\).  Starting with any falling factorial of weight smaller than \(d\), repeatedly apply \cref{lem:consec-span}, always on the fixed target layer \(\ul{w}\), until its weight is \(d\).  This proves the containment from left to right.  The reverse containment is immediate because each displayed falling factorial has degree \(d\).  The rank identity follows from the definition of \(\Ev_{d,w}\).
\end{proof}

\begin{corollary}[Evaluation factorization]\label{cor:evaluation-factorization}
    Let \(d,w\in[0,N]\) with \(d<w\).  Then
    \[
    \Ev_{d,w}=
    \frac{1}{(w-d)!}
    U_{d,d+1}U_{d+1,d+2}\cdots U_{w-1,w}\diagonal_w.
    \]
    In particular, if \(w\le\lfloor N/2\rfloor\), then \(\rank(\Ev_{d,w})=\sqbinom{G}{d}\).

    More generally, if \(d<d'<w\), then
    \[
    \Ev_{d,w}=
    \frac{(w-d')!}{(w-d)!}
    U_{d,d+1}\cdots U_{d'-1,d'}\Ev_{d',w}.
    \]
\end{corollary}

\begin{proof}
    The preceding lemma gives, for every \(r\in[d,w-1]\),
    \[
    \Ev_{r,w}=\frac{1}{w-r}U_{r,r+1}\Ev_{r+1,w}.
    \]
    Iterating this identity and using \(\Ev_{w,w}=\diagonal_w\) gives the first display and the final display.  If \(w\le\lfloor N/2\rfloor\), then \cref{prop:up-matrix}(b) shows that the product of the up matrices has row rank \(\sqbinom{G}{d}\).  The diagonal matrix is invertible and the scalar is nonzero, proving the rank assertion.
\end{proof}

The calculation above applies first to layers in the lower half of the grid.  Reflection transports it to the remaining layers.  The map \(\rho(x)=\top-x\) is a bijection from \(\ul{w}\) to \(\ul{N-w}\).  If a function on \(\ul{w}\) is represented by a polynomial \(P(\mb{X})\) of degree at most \(d\), then the reflected function on \(\ul{N-w}\) is represented by \(P(k_1-1-X_1,\ldots,k_n-1-X_n)\), which has the same degree.  Since reflection is an involution, it follows that \(\Hilbert_d(w)=\Hilbert_d(N-w)\).

\subsection{Proof of the single-layer theorem}

We now combine invertibility on the target layer, the up operator factorization, and reflection.

\begin{proof}[Proof of \cref{thm:single-layer-Hilbert}]
    We consider three cases.

    First suppose \(d\ge w\).  The cumulative evaluation matrix \(\Ev_{[0,d],w}\) contains the invertible block \(\Ev_{w,w}=\diagonal_w\). Hence it has full column rank, and therefore \(\Hilbert_d(w)=\sqbinom{G}{w}\).

    Next suppose \(d<w\le\left\lfloor\frac N2\right\rfloor\). By \cref{cor:full-span,cor:evaluation-factorization}, \(\Hilbert_d(w)=\rank(\Ev_{d,w})=\sqbinom{G}{d}\). By the lower-half unimodality in \cref{prop:up-matrix}(c), \(\sqbinom{G}{d}\le\sqbinom{G}{w}\), and so this is the required minimum.

    Finally suppose \(d<w\) and \(w>\lfloor N/2\rfloor\). Put \(u=N-w\).  Then \(u\le\lfloor N/2\rfloor\), and by reflection and layer symmetry, \(\Hilbert_d(w)=\Hilbert_d(u)\) and \(\sqbinom{G}{w}=\sqbinom{G}{u}\). If \(d\ge u\), then the first case, applied to the layer \(\ul{u}\), gives \(\Hilbert_d(w)=\sqbinom{G}{u}=\sqbinom{G}{w}\). Moreover \(u\le d<N-u\).  If \(d\le\lfloor N/2\rfloor\), lower-half monotonicity gives \(\sqbinom{G}{d}\ge\sqbinom{G}{u}\). If \(d>\lfloor N/2\rfloor\), then \(u<N-d\le\lfloor N/2\rfloor\), and lower-half monotonicity together with symmetry gives \(\sqbinom{G}{d} = \sqbinom{G}{N-d} \ge \sqbinom{G}{u}\). Thus in either case the required minimum is \(\sqbinom{G}{w}\).

    If \(d<u\), then the second case, applied to the layer \(\ul{u}\), gives \(\Hilbert_d(w)=\Hilbert_d(u)=\sqbinom{G}{d}\). Again by lower-half unimodality, \(\sqbinom{G}{d}\le\sqbinom{G}{u}=\sqbinom{G}{w}\), so the required minimum is \(\sqbinom{G}{d}\).  This completes the proof.
\end{proof}

\paragraph*{Why the two branches are necessary.} For the Boolean four-cube \(G=\{0,1\}^4\), the layer sizes are \(1,4,6,4,1\). Taking \(d=4\) and \(w=2\), every function on \(\ul{2}\) is represented by a polynomial of degree at most \(4\), and hence \(\Hilbert_4(2)=6\). The expression \(\min\{\sqbinom{G}{4},\sqbinom{G}{2}\}\) would instead give \(1\).  Thus the branch \(d\ge w\) cannot be omitted.

\paragraph*{Boundary cases.} The proof includes \(d=0\), the singleton layers \(w=0,N\), the boundary \(d=w\), the complete interpolation degree \(d=N\), both parities at the center, one-dimensional grids, and grids with a nontrivial maximum plateau.  No strict-unimodality hypothesis is used.

This determines the Hilbert function of each individual layer.  The closure criterion, however, requires the Hilbert function of arbitrary unions of layers and the interaction among their evaluation blocks.  We turn to that problem next.

	\section{Evaluation transfer and the strong upper-layer rank theorem}\label{sec:upper-layer-transfer}\label{subsec:upper-layer-transfer}

The single-layer calculation determines the rank of each block separately.  For a union of layers, we also need to compare the column spaces of different target-layer blocks at the same row degree.  The following transfer identity gives this comparison directly and proves the second rank input before it is used in the general Hilbert-function argument.

\begin{proposition}[Adjacent evaluation transfer]\label{prop:upper-layer-transfer}
For every \(0\le d\le w<N\),
\begin{equation}\label{eq:adjacent-evaluation-transfer}
\Ev_{d,w+1}
=
\frac{1}{w+1-d}
\Ev_{d,w}
\diagonal_w^{-1}
U_{w,w+1}
\diagonal_{w+1}.
\end{equation}
\end{proposition}

\begin{proof}
We verify the identity entry by entry.  Fix \(\alpha\in\ul d\) and \(\beta\in\ul{w+1}\).  The \((\alpha,\beta)\)-entry of the right-hand side of~\eqref{eq:adjacent-evaluation-transfer} is
\begin{equation}\label{eq:transfer-entry-sum}
\frac{1}{w+1-d}
\sum_{\substack{i\in[n]\\\beta_i>0}}
\beta_i\,
\mb{Y}^{(\alpha)}(\beta-e_i).
\end{equation}
Indeed, the weight-\(w\) predecessors of \(\beta\) are precisely the points \(\beta-e_i\) for which \(\beta_i>0\), and
\[
\frac{\beta!}{(\beta-e_i)!}=\beta_i.
\]

If \(\alpha\nleq\beta\), then \(\alpha\nleq\beta-e_i\) for every predecessor, and every term in~\eqref{eq:transfer-entry-sum} is zero.  Suppose that \(\alpha\le\beta\).  For each \(i\), with both sides interpreted as zero when \(\beta_i=\alpha_i\),
\[
\beta_i\mb{Y}^{(\alpha)}(\beta-e_i)
=
(\beta_i-\alpha_i)\mb{Y}^{(\alpha)}(\beta).
\]
Therefore
\[
\sum_{\substack{i\in[n]\\\beta_i>0}}
\beta_i\mb{Y}^{(\alpha)}(\beta-e_i)
=
\left(
\sum_{i=1}^n(\beta_i-\alpha_i)
\right)
\mb{Y}^{(\alpha)}(\beta)
=
(w+1-d)\mb{Y}^{(\alpha)}(\beta).
\]
After division by \(w+1-d\), the entry is exactly \(\Ev_{d,w+1}(\alpha,\beta)\).
\end{proof}

For \(0\le d\le r<N\), put \(A_{d,r} = \frac{1}{r+1-d} \diagonal_r^{-1} U_{r,r+1} \diagonal_{r+1}\). Then~\cref{prop:upper-layer-transfer} says \(\Ev_{d,r+1}=\Ev_{d,r}A_{d,r}\). Ordinary matrix multiplication and induction therefore give, for every \(d\le w\le v\le N\),
\begin{equation}\label{eq:iterated-evaluation-transfer}
\Ev_{d,v}
=
\Ev_{d,w}
A_{d,w}A_{d,w+1}\cdots A_{d,v-1},
\end{equation}
where the product is empty when \(v=w\).  In particular,
\begin{equation}\label{eq:upper-layer-column-nesting}
\operatorname{col}(\Ev_{d,v})
\subseteq
\operatorname{col}(\Ev_{d,w}).
\end{equation}

\restatethmnow{thm:upper-layer-rank}{Strong upper-layer rank theorem}

Since the block \(\Ev_{d,\min E}\) occurs among the column blocks of \(\Ev_{d,E}\), the equality in \cref{thm:upper-layer-rank} is equivalent to
\begin{equation}\label{eq:upper-layer-column-equality}
\operatorname{col}(\Ev_{d,E})
=
\operatorname{col}(\Ev_{d,\min E}).
\end{equation}

\begin{proof}[Proof of~\cref{thm:upper-layer-rank}]
Put \(w=\min E\).  For every \(v\in E\), we have \(v\ge w\), so~\eqref{eq:upper-layer-column-nesting} gives
\[
\operatorname{col}(\Ev_{d,v})
\subseteq
\operatorname{col}(\Ev_{d,w}).
\]
Every column in the concatenated matrix \(\Ev_{d,E}\) therefore belongs to the column space of its minimum-layer block.  Hence
\[
\operatorname{rank}(\Ev_{d,E})
\le
\operatorname{rank}(\Ev_{d,w}).
\]
Conversely, \(\Ev_{d,w}\) is literally one of the column blocks of \(\Ev_{d,E}\), so
\[
\operatorname{rank}(\Ev_{d,w})
\le
\operatorname{rank}(\Ev_{d,E}).
\]
The two ranks are equal.  Since \(w=\min E\), this proves the theorem.
\end{proof}

The transfer theorem will be used in the next section to split the evaluation matrix of an interval-compatible set into independent rank contributions.

\section{Affine Hilbert functions of general weight-determined sets}\label{sec:Hilbert-general}

By \cref{prop:Hilbert-closure-criterion}, computing finite-degree \(\Z^*\)-closure requires the affine Hilbert function not only of one layer, but of every union of layers.  We now establish this Hilbert function characterization.  Interval-compatible sets isolate the main block-rank calculation, and the Bernasconi-Egidi enumeration then reduces an arbitrary weight set to that calculation.

The strong upper-layer rank theorem was proved in \cref{sec:upper-layer-transfer}.  We now use its column-space form in the interval-compatible rank decomposition.

\subsection{Interval-compatible weight sets}\label{subsec:interval-compatible}

Let \(I=[c,d]\subseteq[0,N]\) be a nonempty interval.  A set of distinct weights \(E_0=\{w_t:t\in I\}\subseteq[0,N]\) is said to be \tsf{\(I\)-compatible} if the following conditions hold:
\begin{enumerate}[{\normalfont(a)}]
    \item \(w_t\ge t\), for every \(t\in I\);
    \item if \(w_t\ne t\), then \(w_t>d\);
    \item if \(s<t\) and \(w_s,w_t\notin I\), then \(w_t<w_s\).
\end{enumerate}
Thus a weight indexed by \(t\) is either fixed at \(t\), or lies above the whole interval; the off-interval weights occur in reverse order.

The following formula is the local rank statement that will later be applied to the paired part of the Bernasconi-Egidi enumeration.

\begin{theorem}[Affine Hilbert functions of interval-compatible weight sets]\label{thm:Hilbert-I-compatible}
Let \(G\) be a uniform grid, let \(I=[c,d]\subseteq[0,N]\), and let \(E_0=\{w_t:t\in I\}\) be \(I\)-compatible.  Then
\[
\Hilbert_d(E_0)
=
\sum_{t\in I}
\min\left\{
\sqbinom{G}{t},
\sqbinom{G}{w_t}
\right\}.
\]
\end{theorem}

The first ingredient discards the evaluation rows below \(I\), leaving the consecutive row interval on which the block-rank calculation is performed.

For \(a,b\in\mb N\), put \(Y_i^{(a\mid b)}=\prod_{j=0}^{b-1}(X_i-(a+j))\), and, for \(\alpha,\beta\in\mb N^n\), put \(\mb Y^{(\alpha\mid\beta)} = \prod_{i=1}^nY_i^{(\alpha_i\mid\beta_i)}\). As usual, a product indexed by the empty set is \(1\).  Then \(\mb Y^{(\alpha)}\mb Y^{(\alpha\mid\beta)} = \mb Y^{(\alpha+\beta)}\).

\begin{lemma}[Interval spanning]\label{lem:interval-spanning}
Let \(E\subseteq[0,N]\) be nonempty, and let \(q\in[0,N]\) satisfy \(q<\min E\).  Then, as spaces of functions on \(\ul E\),
\[
V_q(E)
\subseteq
\spn\left\{
\mb Y^{(\beta)}:
\beta\in\ul{[q+1,q+|E|]}
\right\}.
\]
\end{lemma}

\begin{proof}
We first prove a one-step statement.  Let \(s<\min E\) and \(\alpha\in\ul s\).  Define
\[
P_E(\mb X) = \prod_{w\in E} \left(\sum_{i=1}^nX_i-w\right),
\]
and consider the translated subgrid \(G_\alpha=[\alpha_1,k_1-1]\times\cdots\times[\alpha_n,k_n-1]\). By \cref{prop:filtered-triangular-basis}, the restriction of \(P_E\) to \(G_\alpha\) has the unique expansion
\[
P_E(\mb X)
=
\sum_{\substack{\alpha+\beta\in G\\
\wgt(\beta)\le |E|}}
 c_\beta\mb Y^{(\alpha\mid\beta)}
\,\text{as functions on }G_\alpha.
\]
The degree bound in the summation follows from the filtered part of that proposition.  Evaluating at \(\alpha\), every term with \(\beta\ne0^n\) vanishes, and hence \(c_{0^n}=P_E(\alpha) = \prod_{w\in E}(s-w)\ne0\).

We now multiply the displayed expansion by \(\mb Y^{(\alpha)}\).  If \(x\in\ul E\) and \(x\not\ge\alpha\), then \(\mb Y^{(\alpha)}(x)=0\), and every term \(\mb Y^{(\alpha+\beta)}(x)\) also vanishes because it contains \(\mb Y^{(\alpha)}\) as a factor.  If \(x\ge\alpha\), then \(x\in G_\alpha\) and \(P_E(x)=0\), because \(\wgt(x)\in E\).  Thus the multiplied expansion is valid on all of \(\ul E\), and
\[
0
=
\mb Y^{(\alpha)}P_E
=
\sum_{\substack{\alpha+\beta\in G\\
\wgt(\beta)\le |E|}}
 c_\beta\mb Y^{(\alpha+\beta)}.
\]
Separating the nonzero term with \(\beta=0^n\), we obtain
\[
\mb Y^{(\alpha)}
=
-\frac{1}{c_{0^n}}
\sum_{\substack{\alpha+\beta\in G\\
\beta\ne0^n,\ \wgt(\beta)\le |E|}}
 c_\beta\mb Y^{(\alpha+\beta)}.
\]
Every falling factorial on the right has weight in \([s+1,s+|E|]\).  Therefore
\[
\mb Y^{(\alpha)}
\in
\spn\left\{
\mb Y^{(\gamma)}:
\gamma\in\ul{[s+1,s+|E|]}
\right\}
\quad\text{on }\ul E.
\tag{5.1}
\]

Let now \(\alpha\in G\) have weight \(a\le q\).  Starting with \((5.1)\), replace only the terms in the bottom layer of the current band.  At every step that bottom weight is at most \(q<\min E\), so the one-step statement applies again.  The terms already above the bottom layer remain in the new band.  Repeating this operation \(q-a\) times gives
\[
\mb Y^{(\alpha)}
\in
\spn\left\{
\mb Y^{(\beta)}:
\beta\in\ul{[q+1,q+|E|]}
\right\}.
\]
Since \(E\subseteq[q+1,N]\), we have \(|E|\le N-q\), so the terminal band \([q+1,q+|E|]\) lies in \([0,N]\).  Finally, the restrictions of the falling factorials of weights at most \(q\) span \(V_q(E)\), by \cref{prop:filtered-triangular-basis}.  This proves the lemma.
\end{proof}

The second input controls all target layers lying at or above the row degree.  Its full range is essential in the general decomposition.  By \cref{thm:upper-layer-rank}, equivalently by~\eqref{eq:upper-layer-column-equality}, every higher target-layer block lies in the column space of the minimum-layer block.  The remaining argument is block elimination.  We begin with the elementary rank-splitting fact used in that induction; the proof is immediate.
\begin{lemma}[Block rank splitting]\label{lem:block-rank-splitting}
Let
\[
M=
\begin{bmatrix}
A&B\\
0&F
\end{bmatrix}.
\]
If \(\operatorname{col}(B)\subseteq\operatorname{col}(A)\), then \(\rank M=\rank A+\rank F\).
\end{lemma}
%
%\begin{proof}
%Choose a matrix \(C\) such that \(B=AC\).  Right multiplication by the invertible block matrix
%\[
%\begin{bmatrix}
%I&-C\\
%0&I
%\end{bmatrix}
%\]
%changes \(M\) into
%\[
%\begin{bmatrix}
%A&0\\
%0&F
%\end{bmatrix}.
%\]
%The asserted rank formula follows.
%\end{proof}
The general Bernasconi-Egidi enumeration may leave additional upper layers beyond the interval-compatible core.  The next strengthening shows that those layers are rank-redundant.

\begin{proposition}[Augmented interval-compatible rank]\label{prop:augmented-interval-rank}
Let \(I=[c,d]\subseteq[0,N]\), and let \(E_0=\{w_t:t\in I\}\) be \(I\)-compatible.  Let \(R\subseteq[0,N]\) be disjoint from \(E_0\), and suppose that \(v>\max E_0\) for every \(v\in R\).  Then
\[
\rank(\Ev_{I,E_0\cup R})
=
\sum_{t\in I}
\min\left\{
\sqbinom{G}{t},
\sqbinom{G}{w_t}
\right\}.
\]
In particular, the layers in \(R\) do not increase this rank.
\end{proposition}

\begin{proof}
We induct on \(|I|\).

Suppose first that \(I=\{d\}\).  If \(w_d=d\), then the invertible block \(\Ev_{d,d}=\diagonal_d\) gives \(\Ev_{d,\{d\}\cup R}\) full row rank, so the extra columns indexed by \(R\) cannot increase its rank.  If \(w_d>d\), then \(\{w_d\}\cup R\subseteq[d,N]\) and its minimum is \(w_d\); hence \cref{thm:upper-layer-rank} gives \(\rank(\Ev_{d,\{w_d\}\cup R}) = \rank(\Ev_{d,w_d})\).  In both cases, the result follows from \cref{thm:single-layer-Hilbert}.

Now let \(I=[c,d]\) with \(c<d\), and assume the result for \([c,d-1]\).  Put the column block indexed by \(w_d\) first, and write
\[
\Ev_{[c,d],E_0\cup R}
=
\begin{bmatrix}
\Ev_{d,w_d}&B\\
C&D
\end{bmatrix},
\]
where the lower row block is indexed by \([c,d-1]\).

For \(s=c,c+1,\ldots,d-1\), in this order, replace the row block indexed by \(s\) by that row block minus \(\frac{1}{w_d-s}U_{s,s+1}\) times the row block indexed by \(s+1\).  At the moment the operation indexed by \(s\) is performed, the row block indexed by \(s+1\) has not yet been changed.  Each operation is an invertible elementary block-row operation, and hence preserves rank.  By \cref{cor:evaluation-factorization}, the first column block below the top row becomes zero.  The resulting matrix has the form
\[
\begin{bmatrix}
\Ev_{d,w_d}&B\\
0&F
\end{bmatrix}.
\tag{5.3}
\]

We verify the hypothesis of \cref{lem:block-rank-splitting}.  If \(w_d=d\), then \(\Ev_{d,d}\) is invertible.  Suppose that \(w_d>d\).  Every other core weight is either a fixed weight \(t<d\), whose top block \(\Ev_{d,t}\) is zero, or an off-interval weight larger than \(w_d\), by interval compatibility.  Every element of \(R\) is also larger than \(w_d\).  After the zero blocks are discarded, the target weights in the top row form a subset of \([d,N]\) with minimum \(w_d\).  Equation~\eqref{eq:upper-layer-column-equality} therefore gives
\[
\operatorname{col}(B)
\subseteq
\operatorname{col}(\Ev_{d,w_d}).
\]
Applying \cref{lem:block-rank-splitting} to \((5.3)\), we get
\[
\rank(\Ev_{[c,d],E_0\cup R})
=
\rank(\Ev_{d,w_d})+\rank(F).
\tag{5.4}
\]

It remains to identify \(F\).  For \(s\in[c,d-1]\) and \(v\in(E_0\setminus\{w_d\})\cup R\), the corresponding block is
\[
F_{s,v}
=
\Ev_{s,v}
-
\frac{1}{w_d-s}U_{s,s+1}\Ev_{s+1,v}.
\]
If \(v>s\), then \cref{cor:evaluation-factorization} gives
\[
U_{s,s+1}\Ev_{s+1,v}
=(v-s)\Ev_{s,v}.
\]
If \(v=s\), then \(\Ev_{s+1,s}=0\); if \(v<s\), both evaluation blocks are zero.  Thus, in all cases,
\[
F_{s,v}
=
\frac{w_d-v}{w_d-s}\Ev_{s,v}.
\tag{5.5}
\]
Consequently,
\[
F
=
\operatorname{diag}\left(
\frac{1}{w_d-s}I_{\ul s}:s\in[c,d-1]
\right)
\Ev_{[c,d-1],(E_0\setminus\{w_d\})\cup R}
\operatorname{diag}\left(
(w_d-v)I_{\ul v}:v\ne w_d
\right).
\]
All scalars in the two diagonal factors are nonzero: \(w_d>s\), and the target weights are distinct.  Hence
\[
\rank(F)
=
\rank\left(
\Ev_{[c,d-1],(E_0\setminus\{w_d\})\cup R}
\right).
\tag{5.6}
\]
The set \(E_0\setminus\{w_d\}=\{w_t:t\in[c,d-1]\}\) is \([c,d-1]\)-compatible, and every element of \(R\) remains larger than every weight in this smaller core.  The induction hypothesis applies.  Combining it with \((5.4)\), \((5.6)\), and the single-layer rank formula proves the proposition.
\end{proof}

\begin{proposition}[Interval-compatible rank]\label{prop:interval-compatible-rank}
Let \(I=[c,d]\subseteq[0,N]\), and let \(E_0=\{w_t:t\in I\}\) be \(I\)-compatible.  Then
\[
\rank(\Ev_{I,E_0})
=
\sum_{t\in I}
\min\left\{
\sqbinom{G}{t},
\sqbinom{G}{w_t}
\right\}.
\]
\end{proposition}

\begin{proof}
Take \(R=\varnothing\) in \cref{prop:augmented-interval-rank}.
\end{proof}

\begin{proof}[Proof of \cref{thm:Hilbert-I-compatible}]
The row span of \(\Ev_{I,E_0}\) is contained in \(V_d(E_0)\), because all its row polynomials have degree at most \(d\).  We prove the reverse containment.

If \(c=0\), then \(I=[0,d]\), and the assertion follows from the definition of the degree-\(d\) evaluation space.  Suppose that \(c>0\).  Since \(w_t\ge t\), we have \(c-1<\min E_0\).  Also, \(|E_0|=|I|=d-c+1\). Applying \cref{lem:interval-spanning} with \(q=c-1\), we get
\[
V_{c-1}(E_0)
\subseteq
\spn\left\{
\mb Y^{(\alpha)}:
\alpha\in\ul{[c,d]}
\right\}.
\]
The falling-factorial rows of weights between \(c\) and \(d\) are already indexed by \(I\).  Thus every row of degree at most \(d\) lies in the row span of \(\Ev_{I,E_0}\).  Therefore \(\Hilbert_d(E_0)=\rank(\Ev_{I,E_0})\), and the theorem follows from \cref{prop:interval-compatible-rank}.
\end{proof}

\subsection{The general case}\label{subsec:general-Hilbert}

We now assemble the interval-compatible calculation for an arbitrary weight set.  The \((N,d)\)-BE enumeration records the missing low weights and the present high weights of \(E\subseteq[0,N]\): \([0,d]\setminus E=\{t_\ell<\cdots<t_1\}, \, E\setminus[0,d]=\{w_1<\cdots<w_r\}\).

When some low layers remain unpaired, the next lemma separates their contribution from the terminal interval-compatible part.

\begin{lemma}[Separated lower-layer spanning]\label{lem:separated-lower-spanning}
Let \(I=[c,d]\) with \(1\le c\le d\), let \(E'\subseteq[0,N]\) satisfy \(a+1<c\) for every \(a\in E'\), and let \(E''\) be \(I\)-compatible.  Put \(E=E'\sqcup E''\).  Then
\[
V_d(E)
=
\operatorname{rowspan}(\Ev_{E'\cup I,E}).
\]
\end{lemma}

\begin{proof}
Every row indexed by \(E'\cup I\) has weight at most \(d\), so the right-hand side is contained in \(V_d(E)\).  It remains to show that the rows of weights at most \(c-1\) lie in the displayed row span.  We prove this by induction on \(|E'|\).

If \(E'=\varnothing\), then \(E=E''\).  Since \(E''\) is \(I\)-compatible, \(c-1<\min E''\), and \cref{lem:interval-spanning}, applied with \(q=c-1\), gives
\[
V_{c-1}(E'')
\subseteq
\spn\left\{
\mb Y^{(\alpha)}:
\alpha\in\ul{[c,d]}
\right\},
\]
because \(|E''|=d-c+1\).

Now suppose that \(E'\ne\varnothing\), and let \(a=\min E'\).  Take \(f\in V_{c-1}(E)\).  Since \(a\le c-2\), \cref{thm:single-layer-Hilbert} gives \(V_{c-1}(a)=V(a)\). Since \(\Ev_{a,a}=\diagonal_a\) is invertible, there is a polynomial
\[
P\in
\spn\{\mb Y^{(\alpha)}:\alpha\in\ul a\}
\]
whose restriction to \(\ul a\) is \(f|_{\ul a}\).  On \(E\setminus\{a\}\), the function \(f-P\) still has degree at most \(c-1\).  By induction, it is represented there by a polynomial
\[
Q\in
\spn\left\{
\mb Y^{(\alpha)}:
\alpha\in\ul{(E'\setminus\{a\})\cup I}
\right\}.
\]
Every row weight occurring in \(Q\) is strictly greater than \(a\), and hence every such falling factorial vanishes on \(\ul a\).  Thus \(Q|_{\ul a}=0\).  The polynomial \(P+Q\) represents \(f\) on all of \(E\), and uses only rows indexed by \(E'\cup I\).  This completes the induction.
\end{proof}

We are now ready to prove the Hilbert formula that will be used in \cref{sec:Z-star-clo} through the closure-rank criterion.

\restatethmnow{thm:grid-Hilbert}{Affine Hilbert functions of general weight-determined sets}

\begin{proof}
Let \(m=\min\{\ell,r\}\), and put \(J=(E\cap[0,d])\cup\{t_1,\ldots,t_m\}\subseteq[0,d]\). Let \(I\) be the maximal terminal interval of \([0,d]\) contained in \(J\).  Thus either \(I=[c,d]\) for a uniquely determined \(c\), or \(I=\varnothing\) when \(d\notin J\).  Put \(E'=J\setminus I\).

We record the structure of this decomposition.  Exactly the weights \(t_{m+1},\ldots,t_\ell\) are omitted from \(J\). They are all smaller than every paired missing weight \(t_1,\ldots,t_m\).  Consequently every paired missing weight lies in \(I\).  Moreover, \(E'\subseteq E\cap[0,d]\). If \(I=[c,d]\ne\varnothing\), then
\[
a+1<c\,\text{for every }a\in E'.
\tag{5.7}
\]
If \(r\ge\ell\), then \(m=\ell\), all missing lower weights have been inserted into \(J\), and hence
\[
I=[0,d],
\,
E'=\varnothing.
\tag{5.8}
\]

For each \(t\in I\), define
\[
v_t=
\begin{cases}
 t,&t\in E\cap[0,d],\\
 w_j,&t=t_j\text{ for some }j\in[m].
\end{cases}
\]
Set \(E''=\{v_t:t\in I\}, \, R=\{w_{m+1},\ldots,w_r\}\), where \(R=\varnothing\) when \(m=r\). The second line in the definition of \(v_t\) is \(v_{t_j}=w_j\).  The subscript is the pairing index \(j\), not the missing weight \(t_j\).

If \(I\ne\varnothing\), then \(E''\) is \(I\)-compatible.  Indeed, fixed weights satisfy \(v_t=t\), while every paired weight is above \(d\).  If \(s<t\) and both \(v_s,v_t\notin I\), write \(s=t_j\) and \(t=t_k\).  Then \(s<t\) implies \(j>k\), and so \(v_t=w_k<w_j=v_s\). The correct disjoint decomposition is
\[
E=E'\sqcup E''\sqcup R.
\tag{5.9}
\]
If \(R\ne\varnothing\), then \(r>m\), so \(m=\ell\); in this case every element of \(R\) is larger than every element of \(E''\).

We now distinguish two cases.

\paragraph*{Case 1: \(r\ge\ell\).} By \((5.8)\), we have \(I=[0,d]\), \(E'=\varnothing\), and \(E=E''\sqcup R\). The set \(E''\) is \([0,d]\)-compatible, and all weights in \(R\) lie above the whole core.  Therefore \cref{prop:augmented-interval-rank} gives
\[
\Hilbert_d(E)
=
\rank(\Ev_{[0,d],E''\cup R})
=
\sum_{t=0}^d
\min\left\{
\sqbinom{G}{t},
\sqbinom{G}{v_t}
\right\}.
\]
For \(t\in E\cap[0,d]\), the corresponding term is \(\sqbinom{G}{t}\).  For \(t=t_j\), it is
\[
\min\left\{
\sqbinom{G}{t_j},
\sqbinom{G}{w_j}
\right\}.
\]
Since \(m=\ell=\min\{\ell,r\}\), this is the required formula.  This also includes \(\ell=0\): once every layer in \([0,d]\) is present, all upper layers are redundant at degree \(d\).

\paragraph*{Case 2: \(r<\ell\).} Now \(m=r\), and hence \(R=\varnothing\).  If \(I=\varnothing\), then necessarily \(r=0\), and \(E\subseteq[0,d]\).  Order the row and column blocks of \(\Ev_{E,E}\) by increasing weight.  This matrix is block upper triangular, and its diagonal block at weight \(w\) is the invertible matrix \(\Ev_{w,w}=\diagonal_w\).  Since \(E\subseteq[0,d]\), this invertible matrix occurs among the rows of \(\Ev_{[0,d],E}\).  Hence that matrix has full column rank, and
\[
\Hilbert_d(E)=|\ul E|
=
\sum_{w\in E}\sqbinom{G}{w},
\]
which is the required formula because the second sum is empty.

Assume now that \(I=[c,d]\ne\varnothing\).  If \(c=0\), then \(E'=\varnothing\), and the result follows directly from \cref{thm:Hilbert-I-compatible}.  Suppose that \(c\ge1\).  By \((5.7)\), \cref{lem:separated-lower-spanning} gives
\[
\Hilbert_d(E)
=
\rank(\Ev_{E'\cup I,E'\cup E''}).
\tag{5.10}
\]
Order the row blocks as \(I,E'\), and the column blocks as \(E'',E'\).  Since every row weight in \(I\) is larger than every column-layer weight in \(E'\), \(\Ev_{I,E'}=0\). When the weights in \(E'\) are ordered increasingly, \(\Ev_{E',E'}\) is block upper triangular, with invertible diagonal blocks \(\Ev_{a,a}=\diagonal_a \,(a\in E')\). Thus the matrix in \((5.10)\) has the form
\[
\begin{bmatrix}
\Ev_{I,E''}&0\\
\Ev_{E',E''}&\Ev_{E',E'}
\end{bmatrix}.
\]
Block elimination using the invertible lower-right block gives
\[
\Hilbert_d(E)
=
\rank(\Ev_{I,E''})
+
\sum_{w\in E'}\sqbinom{G}{w}.
\tag{5.11}
\]
By \cref{prop:interval-compatible-rank},
\[
\rank(\Ev_{I,E''})
=
\sum_{t\in I}
\min\left\{
\sqbinom{G}{t},
\sqbinom{G}{v_t}
\right\}.
\]
The fixed terms in \(E'\) and \(I\) account for all weights in \(E\cap[0,d]\), while the paired terms are exactly \((t_j,w_j)\) for \(j\in[r]=[m]\).  Together with \((5.11)\), this gives the required formula.

The two cases exhaust all possibilities, and the proof is complete.
\end{proof}

\paragraph*{Boundary cases.} The proof includes \(E=\varnothing\), \(d=N\), \(\ell=0\), \(r=0\), \(r=\ell\), an empty terminal interval, one-point intervals, excess upper layers, and grids with a nontrivial maximum plateau.  No strict-unimodality hypothesis is used.

The general Hilbert formula is now available for every set appearing in the closure criterion.  To turn the resulting rank-equality test into an explicit closure operator, we must know exactly when the layer-size terms in that formula are equal.  This is the purpose of the next section.

	\section{Layer sizes and maximum plateaux}\label{sec:layer-sequence}

	The Hilbert formula in \cref{thm:grid-Hilbert} reduces closure membership to comparisons among layer sizes.  Symmetry and unimodality are not enough: on a grid with a long central plateau, we must know exactly where equality begins and ends.  We determine that structure in this section by an elementary convolution argument.

	After permuting the coordinates, we may assume that \(k_1\le k_2\le\cdots\le k_n\). For every \(r\in[n]\), let \(G_r=[0,k_1-1]\times\cdots\times[0,k_r-1], \, N_r=\sum_{i=1}^r(k_i-1)\). Put \(N_0=0\).  We use the local abbreviation \(L_r(j)=\sqbinom{G_r}{j}\) and set \(L_r(j)=0\) when \(j\notin[0,N_r]\).  We also put \(L_0(0)=1\) and \(L_0(j)=0\) for \(j\ne0\).  For \(r\in[n]\), define \(p_r=\min\left\{N_{r-1},\left\lfloor\frac{N_r}{2}\right\rfloor\right\}, \, q_r=N_r-p_r\). For the full grid, we write
	\begin{equation}\label{eq:plateau-endpoints}
	p=p_n
	=
	\min\left\{N-\max_{i\in[n]}(k_i-1),\left\lfloor\frac N2\right\rfloor\right\},
	\,
	q=q_n=N-p,
	\end{equation}
	and put \(\lambda=\sqbinom{G}{p}\).

	\subsection{Partial-grid convolution}

	The required layer-size pattern is obtained by adding one coordinate at a time.  For every \(r\in[n]\) and every integer \(j\), separating the last coordinate gives
	\[
	L_r(j)=\sum_{t=0}^{k_r-1}L_{r-1}(j-t).
	\]
	Subtracting the expressions at two consecutive weights gives the useful difference recurrence
	\begin{equation}\label{eq:layer-diff}
	L_r(j+1)-L_r(j)
	=
	L_{r-1}(j+1)-L_{r-1}(j+1-k_r).
	\end{equation}

	\begin{theorem}[Exact layer sequence]\label{thm:layer-sequence}
		For every \(r\in[n]\),
		\[
		L_r(0)<\cdots<L_r(p_r)
		=
		\cdots
		=
		L_r(q_r)>\cdots>L_r(N_r).
		\]
		A strict chain with no adjacent indices is interpreted as empty.  In particular, \([p_r,q_r]\) is the unique interval on which the maximum layer size of \(G_r\) is attained.

		For the full grid, we have
		\[
		\sqbinom{G}{0}<\cdots<\sqbinom{G}{p}
		=
		\cdots
		=
		\sqbinom{G}{q}>\cdots>\sqbinom{G}{N}.
		\]
		If
		\[
		\max_{i\in[n]}(k_i-1)\ge N-\max_{i\in[n]}(k_i-1),
		\]
		then
		\[
		[p,q]
		=
		\left[N-\max_{i\in[n]}(k_i-1),\max_{i\in[n]}(k_i-1)\right],\quad\tx{and }\lambda=\prod_{i=1}^{n-1}k_i,
		\]
		where the coordinates have been ordered as above.
	\end{theorem}
	\begin{proof}
		Coordinatewise reflection is a bijection from the weight-\(j\) layer of \(G_r\) to its weight-\((N_r-j)\) layer.  Hence
		\begin{equation}\label{eq:partial-layer-symmetry}
		L_r(j)=L_r(N_r-j)
		\,(j\in[0,N_r]).
		\end{equation}
		We prove the asserted shape by induction on \(r\).

		For \(r=1\), we have \(L_1(j)=1\) for every \(j\in[0,k_1-1]\).  Since \(p_1=0\) and \(q_1=N_1=k_1-1\), the assertion says precisely that the whole one-dimensional grid is its maximum plateau.

		Now let \(r\ge2\), and assume the result for \(G_{r-1}\).  Put \(A=k_r-1, \, B=N_{r-1}, \, N_r=A+B\). There are two cases.

		\medskip
		\noindent\emph{Case 1: \(A\ge B+1\).} Then \(p_r=B\) and \(q_r=A\).  If \(0\le j<B\), then \(j-A<0\) and \(j+1\in[1,B]\).  Therefore, by~\eqref{eq:layer-diff}, \(L_r(j+1)-L_r(j)=L_{r-1}(j+1)>0\). If \(B\le j<A\), then \(j+1>B\) and \(j-A<0\), so both terms on the right-hand side of~\eqref{eq:layer-diff} vanish.  Hence \(L_r(j+1)-L_r(j)=0\). By~\eqref{eq:partial-layer-symmetry}, the sequence decreases strictly after weight \(A\).  Thus its maximum plateau is exactly \([B,A]\).  For every \(j\in[B,A]\), the convolution window contains every weight of \(G_{r-1}\), and hence
		\[
		L_r(j)=\sum_{s=0}^{B}L_{r-1}(s)=|G_{r-1}|=\prod_{i=1}^{r-1}k_i.
		\]

		\medskip
		\noindent\emph{Case 2: \(A\le B\).} Then
		\[
		p_r=\left\lfloor\frac{N_r}{2}\right\rfloor,
		\,
		q_r=\left\lceil\frac{N_r}{2}\right\rceil.
		\]
		By symmetry, it is enough to prove \(L_r(j+1)-L_r(j)>0 \,(0\le j<p_r)\). If \(j-A<0\), then the second term in~\eqref{eq:layer-diff} is zero.  Also \(j+1\le p_r\le B\), and therefore \(L_r(j+1)-L_r(j)=L_{r-1}(j+1)>0\). Suppose now that \(j-A\ge0\).  Set \(u=j-A, \, v=\min\{j+1,B-j-1\}\). By~\eqref{eq:partial-layer-symmetry} for \(G_{r-1}\), \(L_{r-1}(j+1)=L_{r-1}(v)\). Moreover, \(0\le u<v\le\left\lfloor\frac B2\right\rfloor\). Indeed, if \(j+1\le B/2\), then \(v=j+1\) and \(v-u=A+1>0\).  If \(j+1>B/2\), then \(v=B-j-1\) and \(v-u=N_r-2j-1>0\), because \(j<p_r\).

		The inductive hypothesis shows that \(L_{r-1}\) is nondecreasing on its lower half.  Thus \(L_{r-1}(v)\ge L_{r-1}(u)\). We claim that equality is impossible.  If equality held, the two distinct lower-half weights \(u<v\) would both lie in the maximum plateau of \(G_{r-1}\).  The lower half of that plateau would then contain at least two weights.  Hence \(p_{r-1}<\left\lfloor\frac B2\right\rfloor\), and so, by the definition of \(p_{r-1}\), we would have \(p_{r-1}=N_{r-2}\).  In particular, \(u\ge N_{r-2}\).  Since the side lengths are ordered, \(A=k_r-1\ge k_{r-1}-1\), and therefore \(j=u+A \ge N_{r-2}+k_{r-1}-1 =B\). But \(j<p_r\le B\), a contradiction.  Consequently, \(L_{r-1}(v)>L_{r-1}(u)\), and~\eqref{eq:layer-diff} gives \(L_r(j+1)-L_r(j)>0\).

		This proves strict increase to the center.  Symmetry gives the central equality when \(N_r\) is odd and strict decrease afterward.  The induction is complete.

		It remains only to verify the last assertion about the plateau value when \(A\ge B\).  The case \(A>B\) was already treated in Case~1.  If \(A=B\), then \(p_r=q_r=B\), and the convolution window at weight \(B\) again contains all weights of \(G_{r-1}\).  Thus \(L_r(B)=\prod_{i=1}^{r-1}k_i\). Taking \(r=n\) gives the stated formula for \(G\).
	\end{proof}

	\subsection{Strict unimodality and endpoint maps}

	The convolution theorem identifies the maximum plateau.  We now translate that information into the strict-unimodality criterion and the endpoint maps used by the closure operator.  We call \(G\) strictly unimodal, or \SUt, when the layer sizes increase strictly through the lower central rank; rank symmetry then forces equality of the two central layers when \(N\) is odd.

	\begin{corollary}[\SUt\, criterion]\label{cor:SU2-criterion}
		The following conditions are equivalent:
		\begin{enumerate}[{\normalfont(a)}]
			\item \(G\) is \SUt.
			\item \(q-p\le1\).
			\item \(\max_{i\in[n]}(k_i-1)
			\le
			N-\max_{i\in[n]}(k_i-1)+1\).
			\item \(\max_{i\in[n]}(k_i-1)
			\le
			\left\lceil\frac N2\right\rceil\).
		\end{enumerate}
		After ordering the coordinates, these are also equivalent to \(k_n\le N_{n-1}+2\).
	\end{corollary}
	\begin{proof}
		By~\cref{thm:layer-sequence}, strict unimodality is equivalent to the maximum plateau having the smallest width allowed by rank symmetry, namely \(q-p\le1\).  Put \(A=\max_{i\in[n]}(k_i-1), \, B=N-A\). If \(A\le B\), then \(p=\lfloor N/2\rfloor\) and \(q=\lceil N/2\rceil\), so \(q-p\le1\).  If \(A\ge B+1\), then \(p=B\), \(q=A\), and \(q-p=A-B\).  Hence \(q-p\le1\) if and only if \(A\le B+1\).  This is equivalent to \(2A\le N+1\), and, since \(A\) is an integer, to \(A\le\lceil N/2\rceil\).  The final formulation follows from \(A=k_n-1\) and \(B=N_{n-1}\).
	\end{proof}

	This is precisely the product-of-chains criterion of Dhand quoted in the introduction, translated from chain capacities \(k_i-1\) to the present notation.  For the covering problem, however, the longer plateaux are equally important, and the exact layer-sequence theorem records them as well.

	Recall the endpoint maps \(t^\oplus\) and \(t^\ominus\) defined in the introduction.  The exact layer sequence now gives the form needed in the closure calculation.

	\begin{corollary}[Plateau endpoint maps]\label{cor:oplus-ominus}
		For every \(t\in[0,N]\),
		\[
		t^\ominus=
		\begin{cases}
			p&\text{if }p\le t\le q,\\
			t&\text{otherwise},
		\end{cases}
		\quad\tx{and}\quad
		t^\oplus=
		\begin{cases}
			q&\text{if }p\le t\le q,\\
			t&\text{otherwise}.
		\end{cases}
		\]
		Both maps are order-preserving.  Moreover,
		\[
		(N-t)^\ominus=N-t^\oplus,
		\quad\tx{and}\quad
		(N-t)^\oplus=N-t^\ominus.
		\]
	\end{corollary}
	\begin{proof}
		By~\cref{thm:layer-sequence}, two consecutive layer sizes are equal only when both weights lie in the maximum plateau \([p,q]\).  The displayed formulas follow immediately from the definitions of the endpoint maps.  They show directly that both maps are order-preserving.  The reflection identities follow from \(q=N-p\).
	\end{proof}

	The following comparisons are the precise layer-size input used repeatedly in the closure arguments.

	\begin{corollary}\label{cor:layer-size-comparisons}
		Let \(a,b\in[0,N]\).
		\begin{enumerate}[{\normalfont(a)}]
			\item If \(0\le a\le b\le\left\lfloor\frac N2\right\rfloor\), then \(\sqbinom{G}{a}\le\sqbinom{G}{b}\).
			
			The inequality is strict when \(a<b\le p\).  Equality for distinct lower-half weights occurs exactly when both lie in \([p,\lfloor N/2\rfloor]\).
			\item If \(0\le a\le b\le N\) and \(a+b\le N\), then \(\sqbinom{G}{a}\le\sqbinom{G}{b}\).
			\item For \(1\le i\le p\), we have \(\sqbinom{G}{i}>\sqbinom{G}{i-1}\). For \(p<i\le\left\lfloor\frac N2\right\rfloor\), we have we have \(\sqbinom{G}{i}=\sqbinom{G}{i-1}=\lambda\).
		\end{enumerate}
	\end{corollary}
	\begin{proof}
		Items~(a) and~(c) follow directly from~\cref{thm:layer-sequence}.  For Item~(b), if \(b\le N/2\), use Item~(a).  If \(b>N/2\), then \(a\le N-b<N/2\), and symmetry gives \(\sqbinom{G}{b}=\sqbinom{G}{N-b}\ge\sqbinom{G}{a}\).
	\end{proof}

	\begin{remark}\label{rem:plateau-structure}
		We use the term \tsf{maximum plateau} for the unique interval \([p,q]\) in every uniform grid.  Thus the maximum plateau may be a singleton, the symmetry-forced central pair, or a longer interval.  The grid is \SUt exactly when \(q-p=N\bmod2\). When \(N\) is even, an \SUt grid has a single central maximum.  When \(N\) is odd, it has two equal central maximum layers.  A longer maximum plateau is equivalent to \(q-p\ge2\).
	\end{remark}

	\paragraph*{Boundary examples.} For a one-dimensional grid \([0,k_1-1]\), the entire layer sequence is equal to \(1\).  Consequently, the two-point chain \([0,1]\) is \SUt, while every longer one-dimensional grid is not.

	The grid \([0,2]\times[0,1]\) has layer sequence \(1,2,2,1\). Here \(N=3\), \([p,q]=[1,2]\), and the grid is \SUt.  On the other hand, \([0,3]\times[0,1]\) has layer sequence \(1,2,2,2,1\). Here \(N=4\), \([p,q]=[1,3]\), and the grid is not \SUt.

	The plateau endpoints and layer-size comparisons needed by the Hilbert function criterion are now explicit.  In the next section, we use them to determine \(\Z^*\)-closure and hence the covering thresholds from Proposition~\ref{prop:covering-bridges}.

	\section{Finite-degree \texorpdfstring{\(\Z^*\)-closures}{Z*-closures}}\label{sec:Z-star-clo}

By Proposition~\ref{prop:covering-bridges}, the polynomial covering problem has been reduced to two threshold questions for finite-degree \(\Z^*\)-closure: when the closure first ceases to be all of \([0,N]\), and when it becomes exactly the prescribed weight set.  We now determine the closure itself for every weight-determined set in every uniform grid.  The affine Hilbert formula from \cref{thm:grid-Hilbert} supplies the required rank comparisons, while the exact plateau structure from \cref{thm:layer-sequence} controls the equalities among layer sizes.

Following the problem-to-tool organization established in \cref{sec:prelims}, we first record elementary operator facts, then prove two structural ingredients, combine them in the main closure theorem, and finally turn that theorem into a direct computation.  The earlier characterization recalled in the introduction applies to strictly unimodal uniform grids; the argument below also covers grids with a nontrivial maximum plateau.

Throughout this section, \([p,q]\) denotes the maximum plateau of the layer-size sequence, as in \cref{thm:layer-sequence}; in particular, \(q=N-p\).  We use the endpoint maps \(t^\oplus\) and \(t^\ominus\) from \cref{cor:oplus-ominus}.  For a weight set \(E\subseteq[0,N]\), we write \(N-E=\{N-e:e\in E\}\).

\subsection{Elementary closure facts and the operator \texorpdfstring{\(L_{N,d}\)}{L(N,d)}}

Before proving the two structural ingredients, we record the elementary facts that organize the argument.  Reflection allows the lower and upper ends of the weight interval to be treated together, while the small-set case gives the base of the closure computation.

\begin{proposition}[Reflection equivariance]\label{prop:closure-reflection}
For every \(d\in[0,N]\) and \(E\subseteq[0,N]\),
\[
\zscl_{G,d}(N-E)=N-\zscl_{G,d}(E).
\]
\end{proposition}

\begin{proof}
Let \(\rho(x_1,\ldots,x_n) = (k_1-1-x_1,\ldots,k_n-1-x_n)\). The map \(\rho\) is an affine involution of \(G\), preserves total degree under composition, and maps the layer \(\ul{w}\) onto \(\ul{N-w}\).  Thus a polynomial of degree at most \(d\) vanishes on \(\ul E\) if and only if its pullback by \(\rho\) vanishes on \(\ul{N-E}\).  Consequently, \(\rho\) maps \(\zcl_{G,d}(\ul E)\) onto \(\zcl_{G,d}(\ul{N-E})\).  Taking the complete layers contained in these point closures gives the result.
\end{proof}

\begin{proposition}[Small weight sets]\label{prop:small-zstar}
Let \(d\in[0,N]\).  If \(E\subseteq[0,N]\) and \(|E|\le d\), then
\[
\zscl_{G,d}(E)=E.
\]
\end{proposition}

\begin{proof}
For every \(j\notin E\), consider
\[
P_E(\mb X)
=
\prod_{t\in E}
\left(\sum_{r\in[n]}X_r-t\right).
\]
The empty product is \(1\).  The polynomial \(P_E\) has degree \(|E|\le d\), vanishes on \(\ul E\), and is the nonzero constant \(\prod_{t\in E}(j-t)\) on the layer \(\ul j\).  Hence no weight outside \(E\) belongs to its degree-\(d\) Z*-closure.  Extensivity of closure gives the reverse containment.
\end{proof}

The closure characterization will be expressed by an extensive set operator which records one round of growth forced by the extreme order statistics of the input.  We now recall this operator from the introduction.  For \(d\in[0,N]\) and \(E=\{t_1<\cdots<t_s\}\subseteq[0,N]\), define
\begin{equation}\label{eq:Lnd-definition}
L_{N,d}(E)=
\begin{cases}
E,&s\le d,\\[1mm]
[0,t_{s-d}^{\oplus}]\cup E\cup[t_{d+1}^{\ominus},N],&s\ge d+1.
\end{cases}
\end{equation}
Put \(L_{N,d}^1=L_{N,d}\) and \(L_{N,d}^{r+1}=L_{N,d}\circ L_{N,d}^r \,(r\ge1)\), and define \(\ol L_{N,d}(E)=\bigcup_{r\ge1}L_{N,d}^r(E)\).

\begin{proposition}[Basic properties of \(L_{N,d}\)]\label{prop:Lnd-basic}
Let \(d\in[0,N]\) and \(E,F\subseteq[0,N]\).
\begin{enumerate}[{\normalfont(a)}]
\item \(E\subseteq L_{N,d}(E)\).
\item If \(E\subseteq F\), then \(L_{N,d}(E)\subseteq L_{N,d}(F)\).
\item If \(0\le d<N\), then \(L_{N,d+1}(E)\subseteq L_{N,d}(E)\).
\item \(L_{N,d}(N-E)=N-L_{N,d}(E)\).
\item The increasing chain
\[
E\subseteq L_{N,d}(E)\subseteq L_{N,d}^2(E)\subseteq\cdots
\]
stabilizes.  Its stable value is \(\ol L_{N,d}(E)\), and \(L_{N,d}(\ol L_{N,d}(E))=\ol L_{N,d}(E)\).
Moreover, \(E\mapsto\ol L_{N,d}(E)\) is a closure operator, is monotone decreasing in \(d\), and satisfies
\[
\ol L_{N,d}(N-E)=N-\ol L_{N,d}(E).
\]
\end{enumerate}
\end{proposition}

\begin{proof}
Part~(a) is immediate from \eqref{eq:Lnd-definition}.  For part~(b), the assertion is clear if \(|E|\le d\).  Suppose \(|E|\ge d+1\), and write \(E=\{e_1<\cdots<e_s\}, \, F=\{f_1<\cdots<f_t\}\). The \((d+1)\)-st largest element of \(F\) is at least the \((d+1)\)-st largest element of \(E\), while the \((d+1)\)-st smallest element of \(F\) is at most the corresponding element of \(E\).  Thus \(e_{s-d}\le f_{t-d}, \, f_{d+1}\le e_{d+1}\). The endpoint maps are order-preserving by \cref{cor:oplus-ominus}.  Hence the lower interval added to \(E\) is contained in the lower interval added to \(F\), and the upper interval added to \(E\) is contained in the upper interval added to \(F\).  This proves monotonicity.

For part~(c), only the case \(|E|\ge d+2\) requires attention.  If \(E=\{t_1<\cdots<t_s\}\), then \(t_{s-d-1}\le t_{s-d}, \, t_{d+2}\ge t_{d+1}\). Again the monotonicity of the endpoint maps gives the stated containment.

For part~(d), if \(s\le d\), both sides fix their input.  Otherwise the increasing enumeration of \(N-E\) is \(N-t_s<\cdots<N-t_1\). The two relevant order statistics are \(N-t_{d+1}\) and \(N-t_{s-d}\).  By \cref{cor:oplus-ominus},
\[
(N-t)^\oplus=N-t^\ominus,
\,
(N-t)^\ominus=N-t^\oplus.
\]
Substituting these identities into \eqref{eq:Lnd-definition} proves the reflection formula.

Finally, parts~(a) and~(b) give an increasing chain in the finite set \([0,N]\), so it stabilizes after finitely many strict enlargements.  Its stable value is fixed by \(L_{N,d}\).  Extensivity, monotonicity, and idempotence of the stabilized map follow at once.  Part~(c) may be iterated to obtain degree antitonicity, and part~(d) may be iterated to obtain reflection equivariance of the stabilized map.
\end{proof}

The first structural ingredient will be applied to every \((d+1)\)-subset of the input.  The following proposition shows that the union of the intervals forced by those subsets is exactly one application of \(L_{N,d}\).

\begin{proposition}[Union over \((d+1)\)-subsets]\label{prop:Lnd-union}
Let \(d\in[0,N]\) and \(E\subseteq[0,N]\) with \(|E|\ge d+1\).  Then
\[
L_{N,d}(E)
=
\bigcup_{A\in\binom{E}{d+1}}
\left(
[0,(\min A)^\oplus]\cup A\cup[(\max A)^\ominus,N]
\right).
\]
\end{proposition}

\begin{proof}
Write \(E=\{t_1<\cdots<t_s\}\).  The union of all \((d+1)\)-subsets of \(E\) is \(E\).  The largest possible value of \(\min A\) is \(t_{s-d}\), while the smallest possible value of \(\max A\) is \(t_{d+1}\).  Since the endpoint maps are order-preserving, the union of the lower intervals is \([0,t_{s-d}^\oplus]\), and the union of the upper intervals is \([t_{d+1}^\ominus,N]\).  This is exactly \eqref{eq:Lnd-definition}.
\end{proof}

\begin{lemma}[Full growth]\label{lem:Lnd-full-growth}
Let \(d\in[0,N]\).  If \([0,d]\subseteq F\subseteq[0,N]\), then
\[
L_{N,d}(F)=[0,N].
\]
\end{lemma}

\begin{proof}
The set \(F\) has at least \(d+1\) elements.  Its \((d+1)\)-st smallest element is at most \(d\), and its \(\ominus\)-endpoint is also at most \(d\).  Hence the upper interval inserted by \(L_{N,d}\) begins at or before \(d\).  Together with \([0,d]\subseteq F\), this covers all of \([0,N]\).
\end{proof}

\subsection{The first structural ingredient: the Closure Builder Lemma}

The first structural ingredient says that once \(d+1\) layers are prescribed, the affine Hilbert formula forces complete intervals at both ends of the weight range into the closure.  We first record the layer-size comparison used in its proof.  If \(0\le a\le b\le N-a\), then
\begin{equation}\label{eq:layer-between-reflections}
\sqbinom{G}{a}\le\sqbinom{G}{b}.
\end{equation}
Indeed, \(b\) lies between \(a\) and its reflected weight \(N-a\), and the layer-size sequence cannot fall below its value at either endpoint on this interval.

\begin{lemma}[Closure Builder Lemma]\label{lem:closure-builder}
Let \(d\in[0,N]\) and \(E\subseteq[0,N]\) with \(|E|\ge d+1\).  Then
\[
[0,(\min E)^\oplus]
\cup
[(\max E)^\ominus,N]
\subseteq
\zscl_{G,d}(E).
\]
\end{lemma}

\begin{proof}
We first prove the lower containment.  Choose \(A\subseteq E, \, |A|=d+1, \, \min A=\min E\). By monotonicity of Z*-closure, \(\zscl_{G,d}(A)\subseteq\zscl_{G,d}(E)\). It is therefore enough to prove \([0,(\min A)^\oplus]\subseteq\zscl_{G,d}(A)\). Put \(m=\min A\), and write the \((N,d)\)-BE enumeration as
\[
[0,d]\setminus A=\{t_\ell<\cdots<t_1\},
\,
A\setminus[0,d]=\{w_1<\cdots<w_\ell\}.
\]
The two lists have the same length because \(|A|=d+1\).  By \cref{thm:grid-Hilbert},
\begin{equation}\label{eq:builder-Hilbert}
\Hilbert_d(A)
=
\sum_{u\in A\cap[0,d]}\sqbinom{G}{u}
+
\sum_{j=1}^{\ell}
\min\left\{
\sqbinom{G}{t_j},
\sqbinom{G}{w_j}
\right\}.
\end{equation}
By extensivity, every weight already in \(A\) lies in its closure.  It remains to show that adjoining every missing \(a\le m^\oplus\) leaves this value unchanged.

\paragraph*{Case 1: \(a<m\).} Suppose first that \(a\le d\).  Since no element of \(A\) lies below \(m\), every weight in \([0,m-1]\) is missing.  Thus \(a=t_s, \, s=\ell-a, \, t_j=\ell-j\quad(s\le j\le\ell)\). The weight \(w_j\) is the \(j\)-th member of an \(\ell\)-element subset of \([0,N]\), and hence \(w_j\le N-(\ell-j)\). Therefore \(t_j+w_j\le N\) for \(s\le j\le\ell\).  Since \(t_j\le d<w_j\), \eqref{eq:layer-between-reflections} gives
\[
\sqbinom{G}{t_j}\le\sqbinom{G}{w_j}
\,(s\le j\le\ell).
\]
After adjoining \(a=t_s\), the missing list loses \(t_s\), while the upper list is unchanged and its last member becomes excess.  Applying \cref{thm:grid-Hilbert} before and after the insertion gives
\begin{align*}
\Hilbert_d(A\cup\{a\})-\Hilbert_d(A)
={}&\sqbinom{G}{a}
+
\sum_{j=s}^{\ell-1}
\left(
\sqbinom{G}{t_{j+1}}-\sqbinom{G}{t_j}
\right)
-
\sqbinom{G}{t_\ell}\\
={}&\sqbinom{G}{a}
+
\sqbinom{G}{0}
-
\sqbinom{G}{a}
-
\sqbinom{G}{0}
=0.
\end{align*}
Thus \(a\in\zscl_{G,d}(A)\) by \cref{prop:Hilbert-closure-criterion}.

Suppose now that \(a>d\).  Then every element of \(A\) is larger than \(a\), and hence \(\ell=d+1, \, t_j=d-j+1\). Let \(w_1<\cdots<w_{d+1}\) and \(w'_1<\cdots<w'_{d+2}\) be the upper lists for \(A\) and \(A\cup\{a\}\), respectively.  The elementary order-statistic bounds give \(w_j\le N-(d+1-j), \, w'_j\le N-(d+2-j)\). Thus \(t_j+w_j\le N, \, t_j+w'_j<N\). By \eqref{eq:layer-between-reflections}, every paired minimum in both Hilbert formulas is attained on the \(t_j\)-side.  Therefore
\[
\Hilbert_d(A)
=
\sum_{j=0}^{d}\sqbinom{G}{j}
=
\Hilbert_d(A\cup\{a\}),
\]
and again \(a\) belongs to the closure.

\paragraph*{Case 2: \(m<a\le m^\oplus\).} Such a missing \(a\) can occur only when \(p\le m\le q, \, m^\oplus=q, \, m<a\le q\). This also deals with the range \(m>q\): in that range \(m^\oplus=m\), so there is no additional insertion to prove.

Since \(A\subseteq[m,N]\) and \(|A|=d+1\), \(d\le N-m\le q\). If \(m\le d\), put \(h=\ell-m=|[m,d]\setminus A|\). Then \(t_1,\ldots,t_h\) are precisely the missing weights in \([m,d]\), and all lie on the maximum plateau; moreover, \(t_{h+1}=m-1,\ldots,t_\ell=0\). The first \(h\) upper weights also lie on the plateau.  Indeed, at most \(p\) weights lie above \(q\), and hence \(|A\cap[m,q]| \ge d+1-p \ge d+1-m\). After removing the elements of \(A\cap[m,d]\), this gives \(|A\cap[d+1,q]| \ge h\).

Assume first that \(m<a\le d\).  Then \(a=t_s\) for some \(s\in[h]\).  Removing \(t_s\) from the missing list gives
\[
t'_j=
\begin{cases}
t_j,&j<s,\\
t_{j+1},&s\le j\le\ell-1.
\end{cases}
\]
The upper list is unchanged, but only its first \(\ell-1\) entries remain paired.  Therefore
\begin{align}\label{eq:builder-plateau-difference}
\Hilbert_d(A\cup\{a\})-\Hilbert_d(A)
={}&\sqbinom{G}{a}\notag\\
&+
\sum_{j=s}^{\ell-1}
\left[
\min\left\{\sqbinom{G}{t_{j+1}},\sqbinom{G}{w_j}\right\}
-
\min\left\{\sqbinom{G}{t_j},\sqbinom{G}{w_j}\right\}
\right]\notag\\
&-
\min\left\{\sqbinom{G}{t_\ell},\sqbinom{G}{w_\ell}\right\}.
\end{align}
Suppose first that \(m>0\).  For \(s\le j<h\), all the relevant weights lie on the plateau, so the corresponding terms are zero.  At \(j=h\), the term is \(\sqbinom{G}{m-1}-\lambda\). For \(h+1\le j\le\ell-1\), one has \(t_j=\ell-j\) and \(t_j+w_j\le N\).  Hence both minima in each difference are attained on the \(t\)-side, and the remaining sum telescopes to \(\sqbinom{G}{0}-\sqbinom{G}{m-1} =1-\sqbinom{G}{m-1}\). The final minimum in \eqref{eq:builder-plateau-difference} equals \(\sqbinom{G}{0}=1\), while \(\sqbinom{G}{a}=\lambda\).  Thus the total difference is
\[
\lambda+
\left(\sqbinom{G}{m-1}-\lambda\right)
+
\left(1-\sqbinom{G}{m-1}\right)-1
=0.
\]
If \(m=0\), then \(h=\ell\).  Every \(t_j\), every paired \(w_j\), and the inserted weight \(a\) lies on the plateau.  All shifted summands in \eqref{eq:builder-plateau-difference} vanish, and the added \(+\lambda\) cancels the deleted final pair \(-\lambda\).  Thus the difference is again zero.

It remains to consider \(d<a\le q\).  If \(m>d\), then all weights of \(A\cup\{a\}\) lie above \(d\), and the argument from the second part of Case~1 applies verbatim.  Assume therefore that \(m\le d<a\le q\). Adding \(a\) inserts one plateau weight into the upper list.  For \(1\le j\le h\), both the original and the new \(j\)-th upper weights lie on the plateau, as does \(t_j\); both paired minima are therefore \(\lambda\).  For \(j\ge h+1\), one has \(t_j=\ell-j\).  The original and augmented upper lists satisfy \(w_j\le N-(\ell-j), \, w'_j\le N-(\ell+1-j)\). Hence \(t_j+w_j\le N, \, t_j+w'_j<N\), and \eqref{eq:layer-between-reflections} shows that both paired minima equal \(\sqbinom{G}{t_j}\).  Every summand in the two Hilbert formulas is therefore unchanged.  This proves \([0,m^\oplus]\subseteq\zscl_{G,d}(A) \subseteq\zscl_{G,d}(E)\).

Apply the lower-tail result to the reflected set \(N-E\).  By \cref{prop:closure-reflection,cor:oplus-ominus}, reflection carries the resulting lower interval onto \([(\max E)^\ominus,N]\). This proves the lemma.
\end{proof}

The first structural ingredient now gives one complete direction of the closure characterization.  Applying it to all \((d+1)\)-subsets through \cref{prop:Lnd-union}, and then iterating, yields the forward inclusion.

\begin{corollary}[Forward inclusion]\label{cor:Lbar-in-zstar}
For every \(d\in[0,N]\) and \(E\subseteq[0,N]\),
\[
\ol L_{N,d}(E)\subseteq\zscl_{G,d}(E).
\]
\end{corollary}

\begin{proof}
If \(|E|\le d\), then \(\ol L_{N,d}(E)=E\), and there is nothing to prove.  Suppose \(|E|\ge d+1\).  The Closure Builder Lemma applied to each \((d+1)\)-subset in \cref{prop:Lnd-union}, together with monotonicity of Z*-closure, gives
\[
L_{N,d}(E)\subseteq\zscl_{G,d}(E).
\]
If \(L_{N,d}^r(E)\subseteq\zscl_{G,d}(E)\), the same argument applied to \(L_{N,d}^r(E)\), followed by monotonicity and idempotence of Z*-closure, gives \(L_{N,d}^{r+1}(E)\subseteq\zscl_{G,d}(E)\). The result follows after stabilization.
\end{proof}

\subsection{The second structural ingredient and the closure theorem}

The reverse containment requires low-degree polynomials which keep omitted weights out of the closure.  Symmetric tail sets provide these separating polynomials and form the second structural ingredient.  For \(i\in[0,N]\), define the \tsf{tail set} \(T_{N,i}=[0,i-1]\cup[N-i+1,N]\). Empty intervals are interpreted as empty sets.  Thus \(T_{N,0}=\varnothing\), while \(T_{N,i}=[0,N]\) whenever \(i>\lfloor N/2\rfloor\).

The following Tail Separation Lemma gives exactly the two ranges needed in the closure proof.  Its first part detects when the degree is too small, while its second part supplies the fixed tails used to separate omitted weights.  The range not required in that proof will be determined afterward from the closure theorem.

\begin{lemma}[Tail separation]\label{lem:tail-separation}
Let \(d\in[0,N]\) and \(i\in[0,N]\).
\begin{enumerate}[{\normalfont(a)}]
\item If \(i>d\), then \(\zscl_{G,d}(T_{N,i})=[0,N]\).
\item If \(i\le\min\{d,p\}\), then \(\zscl_{G,d}(T_{N,i})=T_{N,i}\).
\end{enumerate}
\end{lemma}

\begin{proof}
The assertions are immediate when \(i>\lfloor N/2\rfloor\), so assume \(i\le\lfloor N/2\rfloor\).

If \(i>d\), then \([0,d]\subseteq T_{N,i}\).  Hence \cref{thm:grid-Hilbert} gives
\[
\Hilbert_d(T_{N,i})
=
\sum_{s=0}^{d}\sqbinom{G}{s}
=
\Hilbert_d([0,N]).
\]
By \cref{prop:Hilbert-closure-criterion}, the degree-\(d\) Z*-closure is the full weight set.

Now suppose \(i\le\min\{d,p\}\).  Degree antitonicity gives
\[
T_{N,i}
\subseteq
\zscl_{G,d}(T_{N,i})
\subseteq
\zscl_{G,i}(T_{N,i}).
\]
If \(i=0\), then \(T_{N,0}=\varnothing\), and the constant polynomial \(1\) shows that its degree-zero closure is empty.  This proves part~(b) in this case.

Assume \(1\le i\le p\).  At degree \(i\), the only missing low weight is \(i\), and the smallest upper-tail weight is \(N-i+1\).  By symmetry, \(\sqbinom{G}{N-i+1}=\sqbinom{G}{i-1}\). Thus
\[
\Hilbert_i(T_{N,i})
=
\sum_{s=0}^{i-2}\sqbinom{G}{s}
+2\sqbinom{G}{i-1}.
\]
For any \(a\in[i,N-i]\), \cref{thm:layer-sequence} gives \(\sqbinom{G}{a}\ge\sqbinom{G}{i}> \sqbinom{G}{i-1}\). If \(a=i\), all weights in \([0,i]\) are present after adjoining \(a\).  If \(a>i\), the missing low weight \(i\) is paired with \(a\), and the paired minimum is \(\sqbinom{G}{i}\).  In either case,
\[
\Hilbert_i(T_{N,i}\cup\{a\})
=
\sum_{s=0}^{i}\sqbinom{G}{s}
>
\Hilbert_i(T_{N,i}).
\]
By \cref{prop:Hilbert-closure-criterion}, no weight in \([i,N-i]\) lies in the degree-\(i\) Z*-closure.  This proves part~(b).
\end{proof}

The two structural ingredients now meet.  The Closure Builder Lemma controls every interval forced into the closure, while tail separation produces a low-degree polynomial excluding each weight outside the stabilized set.  Together they prove the \(\Z^*\)-closure theorem for every uniform grid stated in the introduction.

\begin{proof}[Proof of~\cref{thm:Zscl-Lbar}]
Fix \(d\in[0,N]\) and \(E\subseteq[0,N]\).  If \(|E|\le d\), then \cref{prop:small-zstar} and \eqref{eq:Lnd-definition} give \(\zscl_{G,d}(E)=E=\ol L_{N,d}(E)\). Assume henceforth that \(|E|\ge d+1\), and put \(F=\ol L_{N,d}(E)\). By \cref{cor:Lbar-in-zstar}, \(F\subseteq\zscl_{G,d}(E)\). It remains to prove the reverse containment.  If \(F=[0,N]\), there is nothing to prove.  Suppose that \(F\) is proper.  By \cref{prop:Lnd-basic}, \(L_{N,d}(F)=F\). Define
\[
i_0=
\max\left\{
 i\in\left[0,\left\lfloor\frac N2\right\rfloor\right]:
 T_{N,i}\subseteq F
\right\}.
\]
This maximum exists because \(T_{N,0}=\varnothing\).  We first note that
\begin{equation}\label{eq:i0-at-most-d}
i_0\le d.
\end{equation}
Indeed, if \(i_0\ge d+1\), then \([0,d]\subseteq T_{N,i_0}\subseteq F\), and \cref{lem:Lnd-full-growth} would give \(L_{N,d}(F)=[0,N]\), a contradiction.

Now fix \(j\notin F\).  Maximality of \(i_0\), together with the properness of \(F\), implies \(T_{N,i_0+1}\nsubseteq F\). Since \(T_{N,i_0+1}\setminus T_{N,i_0}\) consists of the weights \(i_0\) and \(N-i_0\), with repetition removed when they coincide, at least one of these two weights lies outside \(F\).  By \cref{prop:Lnd-basic,prop:closure-reflection}, we may reflect the triple \((E,F,j)\), if necessary, and assume
\begin{equation}\label{eq:i0-outside-F}
i_0\notin F.
\end{equation}
Since \(T_{N,i_0}\subseteq F\), we also have
\begin{equation}\label{eq:j-middle}
j\in[i_0,N-i_0].
\end{equation}
We construct a polynomial of degree at most \(d\) which vanishes on \(\ul E\) but is nonzero at some point of \(\ul j\).

\paragraph*{Case 1: \(i_0\le p\).} By \cref{lem:tail-separation}, \(\zscl_{G,i_0}(T_{N,i_0})=T_{N,i_0}\). Since \(j\notin T_{N,i_0}\), there are a point \(x_0\in\ul j\) and a polynomial \(P\in\mb R[\mb X]\) such that
\[
\deg P\le i_0,
\,
P|_{\ul{T_{N,i_0}}}=0,
\,
P(x_0)\ne0.
\]
We claim that
\begin{equation}\label{eq:counting-case-one}
|E\cap[i_0,N-i_0]|\le d-i_0.
\end{equation}
Suppose otherwise, and let \(A=[N-i_0+1,N]\cup(E\cap[i_0,N-i_0])\). Then \(A\subseteq F\), \(|A|\ge d+1\), and \(\min A\ge i_0\).  Hence \(i_0\in[0,(\min A)^\oplus]\subseteq L_{N,d}(A)\). Monotonicity gives \(L_{N,d}(A)\subseteq L_{N,d}(F)=F\), contradicting \eqref{eq:i0-outside-F}.  This proves \eqref{eq:counting-case-one}.

Define
\[
Q(\mb X)
=
\prod_{t\in E\cap[i_0,N-i_0]}
\left(\sum_{r\in[n]}X_r-t\right).
\]
By \eqref{eq:counting-case-one}, \(\deg(PQ)\le i_0+(d-i_0)=d\). The polynomial \(P\) vanishes on every layer of \(E\) outside \([i_0,N-i_0]\), while \(Q\) vanishes on every layer of \(E\) inside that interval.  Therefore \(PQ|_{\ul E}=0\).  Since \(j\notin F\) and \(E\subseteq F\), one has \(j\notin E\), and so \(Q(x_0)= \prod_{t\in E\cap[i_0,N-i_0]}(j-t)\ne0\). Thus \(PQ(x_0)\ne0\), and \(j\notin\zscl_{G,d}(E)\).

\paragraph*{Case 2: \(i_0>p\).} Since the tails are nested, \(T_{N,p+1}=[0,p]\cup[q,N] \subseteq T_{N,i_0}\subseteq F\). Hence every weight outside \(F\), and in particular \(j\), lies in \([p+1,q-1]\). By \eqref{eq:i0-at-most-d}, we also have \(d\ge p\).  The Tail Separation Lemma at degree \(p\) gives \(\zscl_{G,p}(T_{N,p})=T_{N,p}\). Thus there are \(x_0\in\ul j\) and a polynomial \(P\in\mb R[\mb X]\) such that
\[
\deg P\le p,
\,
P|_{\ul{T_{N,p}}}=0,
\,
P(x_0)\ne0.
\]
We claim that
\begin{equation}\label{eq:counting-case-two}
|E\cap[p,q]|\le d-p.
\end{equation}
If not, let \(A=[q+1,N]\cup(E\cap[p,q])\). Then \(A\subseteq F\), \(|A|\ge d+1\), and \(E\cap[p,q]\ne\varnothing\).  Hence \(\min A\in[p,q]\), so \((\min A)^\oplus=q\).  Therefore \(L_{N,d}(A)\) contains \([0,q]\), while \(A\) contains \([q+1,N]\).  Thus \(L_{N,d}(A)=[0,N]\). Monotonicity and \(L_{N,d}(F)=F\) would imply \(F=[0,N]\), a contradiction.  This proves \eqref{eq:counting-case-two}.

Now put
\[
Q(\mb X)
=
\prod_{t\in E\cap[p,q]}
\left(\sum_{r\in[n]}X_r-t\right).
\]
Then \(\deg(PQ)\le p+(d-p)=d\).  The polynomial \(P\) vanishes on every layer of \(E\) outside \([p,q]\), and \(Q\) vanishes on every layer of \(E\cap[p,q]\).  Hence \(PQ|_{\ul E}=0\).  Since \(j\notin E\),
\[
Q(x_0)=
\prod_{t\in E\cap[p,q]}(j-t)\ne0,
\]
and so \(PQ(x_0)\ne0\).  Thus \(j\notin\zscl_{G,d}(E)\).

In either case, every \(j\notin F\) lies outside \(\zscl_{G,d}(E)\), after transporting the conclusion back through reflection when \eqref{eq:i0-outside-F} was obtained from the reflected triple.  Therefore
\[
\zscl_{G,d}(E)\subseteq F=\ol L_{N,d}(E).
\]
Together with the forward inclusion, this proves the theorem.
\end{proof}

\subsection{Computing the closure: complete tail closures and a direct formula}

The uniform-grid closure theorem identifies \(\Z^*\)-closure with the stabilized set operator.  We now turn that characterization into a direct computation.  We first determine the previously omitted plateau range for every tail set, then characterize the proper fixed points, and finally read off the least fixed superset without iterating \(L_{N,d}\).

\begin{theorem}[Closure of a tail set]\label{thm:tail-closure}
Let \(d\in[0,N]\) and \(i\in[0,N]\).  If \(i>\lfloor N/2\rfloor\), then \(T_{N,i}=[0,N]\).  Assume \(0\le i\le\lfloor N/2\rfloor\).  Then
\[
\zscl_{G,d}(T_{N,i})
=
\begin{cases}
[0,N],&d<i,\\[1mm]
[0,N],&i\le d\le2i-p-1,\\[1mm]
T_{N,i},&d\ge\max\{i,2i-p\}.
\end{cases}
\]
The middle range is empty when its upper endpoint is smaller than its lower endpoint.
\end{theorem}

\begin{proof}
By \cref{thm:Zscl-Lbar}, it is enough to compute \(\ol L_{N,d}(T_{N,i})\).  If \(d<i\), then \([0,d]\subseteq T_{N,i}\), so \cref{lem:Lnd-full-growth} gives the full set in one step.

Assume \(d\ge i\).  If \(2i\le d\), then \(|T_{N,i}|=2i\le d\), and the operator fixes the tail by definition.  It remains to consider \(i\le d\le2i-1\). Put \(x=2i-d-1\). Then \(0\le x\le i-1\), and a direct check of the two relevant order statistics gives
\[
L_{N,d}(T_{N,i})
=
[0,x^\oplus]\cup T_{N,i}\cup[(N-x)^\ominus,N].
\]
If \(x<p\), then \(x^\oplus=x\) and \((N-x)^\ominus=N-x\); both added intervals already lie in the tail, so the tail is fixed.  If \(x\ge p\), then \(x\in[p,q]\), and reflection gives \(N-x\in[p,q]\).  Hence \(x^\oplus=q, \, (N-x)^\ominus=p\), so the one-step image is the full set.  Finally, \(x\ge p \quad\Longleftrightarrow\quad d\le2i-p-1\). This proves all three ranges.
\end{proof}

The tail calculation isolates the degree threshold for the basic symmetric sets.  To compute the closure of an arbitrary input, we next characterize all proper fixed points of the one-step operator.

\begin{proposition}[Proper fixed points]\label{prop:Lnd-fixed-points}
Let \(d\in[0,N]\) and \(F\subsetneq[0,N]\).  Then \(L_{N,d}(F)=F\) if and only if there exists \(i\in[0,\min\{d,p\}]\) such that
\[
T_{N,i}\subseteq F
\,\text{and}\,
|F\setminus T_{N,i}|\le d-i.
\]
\end{proposition}

\begin{proof}
Suppose first that \(L_{N,d}(F)=F\).  If \(|F|\le d\), choose \(i=0\).  Otherwise write \(F=\{t_1<\cdots<t_s\}, \, k=s-d\). The lower interval inserted by \(L_{N,d}\) contains \([0,t_k]\), and \(t_k\ge k-1\).  Since \(F\) is fixed, \([0,k-1]\subseteq F\). Similarly, the \(k\)-th largest element is \(t_{d+1}\le N-k+1\), and the upper interval inserted by the operator contains \([N-k+1,N]\).  Hence \(T_{N,k}\subseteq F\). We claim that \(k\le p\).  If \(k\ge p+1\), then \(t_k\ge k-1\ge p\), so the inserted lower interval reaches at least \(q\).  Also \(t_{d+1}\le N-k+1\le q\), so the inserted upper interval begins at or before \(p\).  The two intervals would cover \([0,N]\), forcing \(F=[0,N]\), a contradiction.  Thus \(k\le p\), and \(|F\setminus T_{N,k}| =s-2k =d-k\). Since \(k=s-d\le d\) whenever this displayed quantity is nonnegative, \(k\in[0,\min\{d,p\}]\).  More directly, the equality \(d-k=|F\setminus T_{N,k}|\ge0\) gives \(k\le d\).

Conversely, suppose \(i\le\min\{d,p\}\), \(T_{N,i}\subseteq F\), and \(|F\setminus T_{N,i}|\le d-i\). If \(|F|\le d\), the operator fixes \(F\) by definition.  Otherwise put \(k=|F|-d\).  Since \(|T_{N,i}|=2i\), \(|F| =2i+|F\setminus T_{N,i}| \le d+i\), so \(1\le k\le i\le p\).  Hence \(T_{N,k}\subseteq F\).  The \(k\)-th smallest element of \(F\) is \(k-1\), and its \(k\)-th largest element is \(N-k+1\).  Since \(k\le p\), neither endpoint is changed by \(\oplus\) or \(\ominus\).  Thus \(L_{N,d}\) adds only \(T_{N,k}\), which is already contained in \(F\).  Therefore \(L_{N,d}(F)=F\).
\end{proof}

Once the proper fixed points are known, the least fixed superset of the input can be found by a single scan through the possible tails.  This removes the need to iterate the operator when computing a closure.

\begin{theorem}[Direct closure formula]\label{thm:closure-direct}
Let \(E\subseteq[0,N]\) and \(d\in[0,N]\).  Examine the integers \(0\le i\le\min\{d,p\}\) in increasing order.  If there is an \(i\) satisfying
\begin{equation}\label{eq:direct-witness}
E\cup T_{N,i}\ne[0,N]
\,\text{and}\,
|E\setminus T_{N,i}|\le d-i,
\end{equation}
let \(i_*\) be the least such integer.  Then \(\zscl_{G,d}(E)=E\cup T_{N,i_*}\). If no such integer exists, then \(\zscl_{G,d}(E)=[0,N]\).
\end{theorem}

\begin{proof}
By \cref{thm:Zscl-Lbar}, it is enough to compute the least fixed point of \(L_{N,d}\) containing \(E\).  Such a least fixed point is \(\ol L_{N,d}(E)\): every fixed superset of \(E\) contains every iterate by monotonicity.

Suppose \(i\) satisfies \eqref{eq:direct-witness}, and put \(F_i=E\cup T_{N,i}\). The set \(F_i\) is proper, contains \(T_{N,i}\), and satisfies \(|F_i\setminus T_{N,i}|=|E\setminus T_{N,i}|\le d-i\). Thus \cref{prop:Lnd-fixed-points} shows that \(F_i\) is fixed.

Assume that at least one witness exists, and let \(i_*\) be the least one.  Let \(H\) be any fixed superset of \(E\).  If \(H=[0,N]\), then it contains \(E\cup T_{N,i_*}\).  If \(H\) is proper, \cref{prop:Lnd-fixed-points} supplies \(j\in[0,\min\{d,p\}]\) with \(T_{N,j}\subseteq H, \, |H\setminus T_{N,j}|\le d-j\). Then \(E\cup T_{N,j}\subseteq H\ne[0,N]\) and \(|E\setminus T_{N,j}| \le |H\setminus T_{N,j}| \le d-j\). Thus \(j\) is also a witness, and so \(j\ge i_*\).  Since the tails are nested, \(T_{N,i_*}\subseteq T_{N,j}\subseteq H\). Therefore every fixed superset of \(E\) contains \(E\cup T_{N,i_*}\), while that set is itself fixed.  It is the least fixed superset.

If no witness exists and \(\ol L_{N,d}(E)\) were proper, the fixed-point criterion applied to that proper fixed set would produce a witness, a contradiction.  Hence the stable value, and therefore the Z*-closure, is the full weight set.
\end{proof}

\subsection{Consequences of the direct formula}

The direct formula makes the closure operator completely explicit and also shows that it remembers considerably less about the grid than the affine Hilbert function does.

\begin{corollary}[Closure-profile invariance]\label{cor:profile-invariance}
Let \(G\) and \(G'\) be uniform grids with the same total capacity \(N\) and the same maximum-plateau start \(p\).  Then, for every \(d\in[0,N]\) and every \(E\subseteq[0,N]\),
\[
\zscl_{G,d}(E)=\zscl_{G',d}(E).
\]
\end{corollary}

\begin{proof}
The criterion in \cref{thm:closure-direct} uses the grid only through \(N\), the degree \(d\), and the allowed range \(0\le i\le p\).  The tails \(T_{N,i}\) are the same for the two grids, so the least witness, or its nonexistence, is the same.
\end{proof}

\begin{remark}
The affine Hilbert functions of the two grids need not agree.  The formula in \cref{thm:grid-Hilbert} depends on the complete sequence \(\sqbinom{G}{0},\ldots,\sqbinom{G}{N}\), not only on \(N\) and \(p\).
\end{remark}

\begin{corollary}[Monotonicity with respect to plateau width]\label{cor:plateau-monotonicity}
Fix \(N\).  Let \(G_1\) and \(G_2\) be uniform grids with maximum plateaux \([p_1,N-p_1]\) and \([p_2,N-p_2]\), respectively, where
\[
0\le p_1\le p_2\le\left\lfloor\frac N2\right\rfloor.
\]
Then, for every \(d\in[0,N]\) and every \(E\subseteq[0,N]\),
\[
\zscl_{G_2,d}(E)
\subseteq
\zscl_{G_1,d}(E).
\]
Thus a wider maximum plateau gives a weakly larger finite-degree Z*-closure.
\end{corollary}

\begin{proof}
Every witness index available for \(p_1\) is also available for \(p_2\).  If there is no witness for \(p_1\), the closure for \(G_1\) is the full set, and the containment is immediate.  Otherwise let \(i_1\) be the least witness for \(p_1\).  It is also a witness for \(p_2\), so the least witness \(i_2\) for \(p_2\) satisfies \(i_2\le i_1\).  Since \(T_{N,i_2}\subseteq T_{N,i_1}\), \cref{thm:closure-direct} gives \(E\cup T_{N,i_2} \subseteq E\cup T_{N,i_1}\). This is the desired containment.
\end{proof}

The closure operator is now fully characterized and directly computable.  By Proposition~\ref{prop:covering-bridges}, its two threshold events are exactly the two polynomial-covering degrees.  We therefore return to the original problem.

\section{Solving the polynomial covering problems}\label{sec:covering}

We now return to Problem~\ref{prob:polynomial-covering}.  Proposition~\ref{prop:covering-bridges} showed that \(\PC_G(E)\) is the first degree at which the \(\Z^*\)-closure is proper, whereas \(\PPC_G(E)\) is the first degree at which the closure equals \(E\).  Since \cref{thm:closure-direct} computes that closure explicitly, the covering problem has been reduced to reading off these two thresholds.  We do this first, then derive certifying degrees and a single linear-time algorithm.

\subsection{Exact covering formulas}

The direct closure characterization gives the first answer to Problem~\ref{prob:polynomial-covering}: a pair of exact covering formulas.  We introduce the admitting-set terminology from the strictly-unimodal covering problem, with the necessary plateau restriction.  Let \(d\in[0,N]\) and \(0\le i\le\min\{d,p\}\). We say that \(E\subseteq[0,N]\) is \tsf{\((d,i)\)-admitting} if
\begin{equation}\label{eq:admitting-definition}
E\cup T_{N,i}\ne[0,N]
\,\text{and}\,
|E\setminus T_{N,i}|\le d-i.
\end{equation}
We say that \(E\) is \tsf{\(d\)-admitting} if it is \((d,i)\)-admitting for at least one \(i\in[0,\min\{d,p\}]\).

The two thresholds from Proposition~\ref{prop:covering-bridges} can now be read directly from \cref{thm:closure-direct,prop:Lnd-fixed-points}.  This gives the two main answers to Problem~\ref{prob:polynomial-covering}.

\begin{theorem}[Nontrivial polynomial covers]\label{thm:PC-formula}
For every \(E\subsetneq[0,N]\),
\begin{equation}\label{eq:PC-admitting-formula}
\PC_G(E)
=
\min\left\{
 d\in[0,N]:E\text{ is }d\text{-admitting}
\right\}.
\end{equation}
Equivalently,
\begin{equation}\label{eq:PC-formula}
\PC_G(E)=
\min_{\substack{0\le i\le p\\ E\cup T_{N,i}\ne[0,N]}}
\left(i+|E\setminus T_{N,i}|\right).
\end{equation}
\end{theorem}

\begin{proof}
By~\cref{prop:covering-bridges}, \(\PC_G(E)\) is the first degree for which the Z*-closure is proper.  By~\cref{thm:closure-direct}, this is precisely the first \(d\) for which \(E\) is \(d\)-admitting.  This proves~\eqref{eq:PC-admitting-formula}.

For a fixed index \(i\in[0,p]\) satisfying \(E\cup T_{N,i}\ne[0,N]\), properness gives \(|E\setminus T_{N,i}|\le N-2i\).  Hence \(d=i+|E\setminus T_{N,i}|\) lies in \([i,N]\) and is the least degree for which~\eqref{eq:admitting-definition} holds.  Taking the minimum over \(0\le i\le p\) proves~\eqref{eq:PC-formula}.  The index \(i=0\) is always eligible because \(E\ne[0,N]\), so the minimum exists.
\end{proof}

\begin{theorem}[Proper polynomial covers]\label{thm:PPC-formula}
For every \(E\subsetneq[0,N]\),
\begin{equation}\label{eq:PPC-formula}
\PPC_G(E)
=
|E|-
\max\{i\in[0,p]:T_{N,i}\subseteq E\}.
\end{equation}
\end{theorem}

\begin{proof}
Since \(L_{N,d}\) is extensive,
\[
\ol L_{N,d}(E)=E
\quad\Longleftrightarrow\quad
L_{N,d}(E)=E.
\]
By~\cref{prop:Lnd-fixed-points}, this happens precisely when there is an index \(i\in[0,p]\) such that
\[
T_{N,i}\subseteq E
\,\text{and}\,
|E\setminus T_{N,i}|\le d-i.
\]
For such an index, \(|T_{N,i}|=2i\), and the least admissible degree is \(i+|E\setminus T_{N,i}| =i+(|E|-2i) =|E|-i\).  Since \(T_{N,i}\subseteq E\), one has \(|E|\ge2i\), so this degree is at least \(i\) and is admissible.  It is minimized by the largest symmetric tail contained in \(E\), proving~\eqref{eq:PPC-formula}.
\end{proof}

The restriction \(i\le p\) in both formulas is essential.  It is the only change needed to pass from the strictly-unimodal formulas to grids with a longer maximum plateau.

\begin{corollary}[Tail covers]\label{cor:tail-cover}
For every \(0\le i\le\lfloor N/2\rfloor\),
\begin{equation}\label{eq:tail-cover}
\PC_G(T_{N,i})
=
\PPC_G(T_{N,i})
=
2i-\min\{i,p\}
=
\max\{i,2i-p\}.
\end{equation}
\end{corollary}

\begin{proof}
The proper-cover formula gives \(\PPC_G(T_{N,i})=2i-\min\{i,p\}\). For the nontrivial-cover formula, consider a witness \(j\in[0,p]\).  If \(j\le i\), then \(j+|T_{N,i}\setminus T_{N,j}| =j+2(i-j) =2i-j\), which is minimized at \(j=\min\{i,p\}\).  If \(j>i\), which can occur only when \(i<p\), then the candidate equals \(j\), and this does not improve on the valid choice \(j=i\).  Hence the two covering degrees have the same value.
\end{proof}

\begin{remark}[Strictly-unimodal specialization]\label{rem:covering-SU2}
If \(G\) is an \SUt\,grid, then \(p=\left\lfloor\frac N2\right\rfloor\). Thus the restriction \(i\le p\) imposes no additional condition on the symmetric tails used in the earlier theory.  The two formulas above reduce exactly to the earlier strictly-unimodal formulas discussed in the introduction.

On a longer plateau, indices larger than \(p\) cannot be used.  For example, on the one-dimensional grid \(G=[0,2]\), one has \(N=2\), \(p=0\), and, for \(E=\{0,2\}\), \(\PC_G(E)=\PPC_G(E)=2\). Allowing the invalid witness \(i=1>p\) would incorrectly give degree \(1\).
\end{remark}

The formulas above solve the polynomial-covering problem for every uniform grid.  They concern polynomial covers; we do not assert an equality with hyperplane-cover parameters, or an exact-cover formula, for a general uniform grid.  The nontrivial-cover threshold also has a second interpretation, which we record next.

\subsection{Certifying degrees}

The same degree which first permits a nontrivial polynomial cover also first permits a nonconstant polynomial function which is constant on the prescribed set.  Let \(S\subsetneq G\).  A polynomial \(P\in\mb R[\mb X]\) is a \tsf{certifying polynomial} for \(S\) if its induced polynomial function on \(G\) is nonconstant, while \(P|_S\) is constant.  The \tsf{certifying degree} of \(S\), denoted by \(\certdeg(S)\), is the least degree of a certifying polynomial for \(S\).

\begin{theorem}[Certifying degree]\label{thm:certdeg}
For every nonempty proper weight set \(E\subsetneq[0,N]\),
\begin{equation}\label{eq:certdeg-PC}
\certdeg(\ul E)=\PC_G(E).
\end{equation}
For the empty set,
\begin{equation}\label{eq:certdeg-empty}
\certdeg(\emptyset)=1,
\,
\PC_G(\emptyset)=\PPC_G(\emptyset)=0.
\end{equation}
\end{theorem}

\begin{proof}
Let \(E\ne\emptyset\).  Suppose that \(P\) is a certifying polynomial of degree \(d\), and that \(P|_{\ul E}=c\).  Then \(P-c\) vanishes on \(\ul E\) and is not the zero function on \(G\), because the function induced by \(P\) is nonconstant.  Thus \(P-c\) is a nontrivial polynomial cover, and \(\PC_G(E)\le\certdeg(\ul E)\). Conversely, a nontrivial cover of the nonempty set \(\ul E\) vanishes there but is not the zero function on \(G\).  Since it takes the value \(0\) somewhere and a nonzero value elsewhere, it is a nonconstant polynomial function, and hence is a certifying polynomial.  This proves~\eqref{eq:certdeg-PC}.

For \(E=\emptyset\), the constant polynomial \(1\) is both a nontrivial and a proper cover, so the two covering degrees are zero.  A certifying polynomial is required to be nonconstant, whereas every degree-zero polynomial function is constant; hence \(\certdeg(\emptyset)\ge1\).  Since every side length is at least two, a coordinate function is nonconstant on \(G\), and therefore \(\certdeg(\emptyset)=1\).
\end{proof}

Combining~\cref{thm:PC-formula,thm:certdeg}, for every nonempty proper \(E\) we have
\[
\certdeg(\ul E)
=
\min_{\substack{0\le i\le p\\E\cup T_{N,i}\ne[0,N]}}
\left(i+|E\setminus T_{N,i}|\right).
\]

\subsection{Direct linear-time algorithms}

The exact formulas solve the numerical covering problem, and the direct closure theorem makes them computable without constructing the intermediate iterates of \(L_{N,d}\).  The final part of the solution is a single linear-time procedure which computes the closure and all three degree parameters.  Suppose that \(E\subsetneq[0,N]\) is given by its indicator vector \((e_0,e_1,\ldots,e_N)\in\{0,1\}^{N+1}\). For \(0\le i\le p\), put
\begin{equation}\label{eq:central-count}
c_i=|E\setminus T_{N,i}|=|E\cap[i,N-i]|.
\end{equation}
These numbers satisfy
\begin{equation}\label{eq:central-count-update}
c_0=|E|,
\,
c_{i+1}=c_i-e_i-e_{N-i}
\quad(0\le i<p).
\end{equation}
The two subtracted positions are distinct because \(i<p\le\lfloor N/2\rfloor\).  Moreover,
\begin{equation}\label{eq:proper-tail-test}
E\cup T_{N,i}\ne[0,N]
\quad\Longleftrightarrow\quad
c_i<N-2i+1.
\end{equation}

The following pseudocode computes the degree-\(d\) \(\Z^*\)-closure, \(\PC_G(E)\), \(\PPC_G(E)\), and the certifying degree in \(O(n+N)\) time.

\begin{algorithm}[H]
\DontPrintSemicolon
\KwIn{The side lengths \(k_1,\ldots,k_n\), a degree \(d\in[0,N]\), and the indicator vector \((e_0,\ldots,e_N)\) of a proper set \(E\subsetneq[0,N]\).}
\KwOut{The indicator vector of \(\zscl_{G,d}(E)\), and the values \(\PC_G(E)\), \(\PPC_G(E)\), and \(\certdeg(\ul E)\).}
Compute \(N=\sum_{r\in[n]}(k_r-1)\) and
\(p=\min\{N-\max_r(k_r-1),\lfloor N/2\rfloor\}\)\;
\(s\leftarrow\sum_{j=0}^Ne_j\), \(c\leftarrow s\), \(i_*\leftarrow\bot\), \(D\leftarrow N+1\)\;
\For{\(i\leftarrow0\) \KwTo \(p\)}{
  \If{\(c<N-2i+1\)}{
    \(D\leftarrow\min\{D,i+c\}\)\;
    \If{\(i_* = \bot\), \(i\le d\), and \(c\le d-i\)}{
      \(i_*\leftarrow i\)\;
    }
  }
  \If{\(i<p\)}{
    \(c\leftarrow c-e_i-e_{N-i}\)\;
  }
}
\If{\(i_*=\bot\)}{
  \(F\leftarrow[0,N]\)\;
}
\Else{
  \(F\leftarrow E\cup T_{N,i_*}\)\;
}
\(\tau\leftarrow0\)\;
\While{\(\tau<p\) and \(e_\tau=e_{N-\tau}=1\)}{
  \(\tau\leftarrow\tau+1\)\;
}
\(D_{\mathrm{proper}}\leftarrow s-\tau\)\;
\eIf{\(s=0\)}{
  \(D_{\mathrm{cert}}\leftarrow1\)\;
}{
  \(D_{\mathrm{cert}}\leftarrow D\)\;
}
\Return the indicator vector of \(F\), and \(D,D_{\mathrm{proper}},D_{\mathrm{cert}}\)\;
\caption{Computing closure, covering degrees, and certifying degree}\label{algo:covering-data}
\end{algorithm}

\begin{theorem}[Correctness and complexity]\label{thm:covering-algorithms}
For every input satisfying the hypotheses of Algorithm~\ref{algo:covering-data}, the algorithm returns
\[
F=\zscl_{G,d}(E),
\,
D=\PC_G(E),
\,
D_{\mathrm{proper}}=\PPC_G(E),
\,
D_{\mathrm{cert}}=\certdeg(\ul E).
\]
Its running time is \(O(n+N)\), and it uses \(O(1)\) auxiliary space apart from the input and output vectors.
\end{theorem}

\begin{proof}
By~\eqref{eq:central-count-update}, the loop maintains \(c=c_i\).  The first test in the loop is exactly the properness condition in~\eqref{eq:proper-tail-test}.  Among the indices satisfying it, the update \(D\leftarrow\min\{D,i+c_i\}\) computes~\cref{thm:PC-formula}.  The first index which also satisfies \(i\le d\) and \(c_i\le d-i\) is the least witness in~\cref{thm:closure-direct}; hence the returned set \(F\) is the degree-\(d\) Z*-closure.

Since the tails are nested, \(T_{N,i}\subseteq E\) holds precisely for \(i=0,1,\ldots,\tau\), where \(\tau\) is the value produced by the inward scan.  Therefore~\cref{thm:PPC-formula} gives \(D_{\mathrm{proper}}=|E|-\tau=\PPC_G(E)\). The rule for \(D_{\mathrm{cert}}\) is exactly~\cref{thm:certdeg}.

Computing \(N\), the largest coordinate capacity, and \(p\) takes \(O(n)\) time.  Each of the two scans uses at most \(p+1\le N+1\) iterations, and writing the output indicator vector takes \(O(N)\) time.  Only a constant number of integer counters is stored besides the input and output.  This proves the stated bounds.
\end{proof}
The algorithm shows that repeatedly applying \(L_{N,d}\) is unnecessary.  The direct closure formula identifies the final closure from the first admissible tail.  The linear scan above therefore computes the stabilized closure without constructing its intermediate iterates.  The same closure scan also applies to \(E=[0,N]\), in which case it returns \([0,N]\); only the three covering and certification outputs require the standing assumption that \(E\) is proper.  Thus the original polynomial-covering problem, its proper-cover variant, and the certification problem are all solved through the closure characterization.

\section{Further remarks}\label{sec:further-remarks}

\paragraph*{A positive-characteristic obstruction.}\label{subsec:positive-char-obstruction}

The characteristic-zero transfer theorem cannot be extended unchanged to positive characteristic.  Consider the Boolean three-cube, the weight set \(E=\{0,2\}\), and degree \(d=1\).  With rows indexed by \(1,X_1,X_2,X_3\) and columns indexed by \(000,\quad110,\quad101,\quad011\), the evaluation matrix is
\[
M=
\begin{pmatrix}
1&1&1&1\\
0&1&1&0\\
0&1&0&1\\
0&0&1&1
\end{pmatrix}.
\]
Its determinant is \(-2\).  Hence \(\operatorname{rank}_{\mb{Q}}(M)=4\). Over \(\mb{F}_2\), the last three rows sum to zero, while the minor formed by the first three rows and first three columns is nonzero.  Therefore \(\operatorname{rank}_{\mb{F}_2}(M)=3\). Thus even the degree-one affine Hilbert function of this layer union depends on the characteristic.  In particular, the characteristic-zero formula of \cref{thm:grid-Hilbert} already fails over \(\mb{F}_2\).  We do not attempt a positive-characteristic classification here.

\paragraph*{Why arbitrary nonuniform grids are different.}\label{subsec:nonuniform-grid}

For arbitrary coordinate sets, we notice that even the degree-one behavior may depend on the actual coordinates and not only on the index weights.
\begin{remark}[A degree-one diagnostic]\label{rem:nonuniform-grid}
Let \(G=\{0,1,2\}^2, \, H=\{0,1,3\}^2\), where the weight on \(H\) is the sum of the coordinate indices.  The two grids have the same layer sizes.  At degree \(d=1\), however, their affine Hilbert functions behave differently.

The weight-two layer of \(G\) is \((2,0),\quad(1,1),\quad(0,2)\). These points lie on the line \(X_1+X_2=2\), and they are not all equal.  Hence \(\Hilbert_1^G(\{2\})=2\). Every point in a layer of weight \(a\ne2\) lies off this line.  Therefore \(\Hilbert_1^G(\{a,2\})=3> \Hilbert_1^G(\{2\}) \,(a\in\{0,1,3,4\})\).

The weight-two layer of \(H\) is \((3,0),\quad(1,1),\quad(0,3)\). These three points are affinely independent, since
\[
\det
\begin{pmatrix}
1&3&0\\
1&1&1\\
1&0&3
\end{pmatrix}
=-3\ne0.
\]
Thus \(\Hilbert_1^H(\{2\})=3\). The whole degree-one function space has dimension three, so adjoining any other layer cannot increase the rank: \(\Hilbert_1^H(\{a,2\}) = \Hilbert_1^H(\{2\}) \,(a\in\{0,1,3,4\})\). The qualifier \(d=1\) is essential: the strict inequality over \(G\) fails at degree zero, and the equality over \(H\) need not hold in higher degree.
\end{remark}
\noindent This example explains why the exact formulas of the paper are stated for uniform grids and their coordinatewise affine arithmetic-progression images, rather than for arbitrary finite Cartesian products.

\paragraph*{Further directions.}\label{subsec:scope-limits}

The following directions have not been considered in this work, and we believe the techniques used in this work will not sufficen, and we need new ideas.
\begin{itemize}
\item a classification in positive characteristic.
\item arbitrary nonuniform coordinate sets.
\item exact polynomial covers.
\item \emph{higher-multiplicity} polynomial covers.
\end{itemize}

\paragraph*{Acknowledgements.}  The author thanks:
	\begin{itemize}
		\item  Srikanth Srinivasan for critical feedback on a preliminary version of this paper.

		\item  Murali K. Srinivasan for some very enlightening discussions on his work~\cite{srinivasan2011symmetric} as well as the works of Canfield~\cite{canfield1980sperner}, Proctor, Saks, and Sturtevant~\cite{proctor1980product}, and Proctor~\cite{proctor1982representations}.

		\item  Lajos R\'onyai for pointers to some relevant literature.
	\end{itemize}

\cleardoublepage
\raggedright
\printbibliography

@article{abbe-shpilka-ye,
  author = {Abbe, Emmanuel and Shpilka, Amir and Ye, Min},
  title = {{Reed--Muller Codes: Theory and Algorithms}},
  journal = {IEEE Transactions on Information Theory},
  year = {2021},
  volume = {67},
  number = {6},
  pages = {3251--3277},
  doi = {10.1109/TIT.2020.3004749},
  addendum = {\url{https://doi.org/10.1109/TIT.2020.3004749}}
}

@article{alon-1999-combinatorial,
  author = {Alon, Noga},
  title = {{Combinatorial Nullstellensatz}},
  journal = {Combinatorics, Probability and Computing},
  year = {1999},
  volume = {8},
  number = {1--2},
  pages = {7--29},
  doi = {10.1017/S0963548398003411},
  addendum = {\url{https://doi.org/10.1017/S0963548398003411}}
}

@book{babai-frankl-linear,
  author = {Babai, L{\'a}szl{\'o} and Frankl, P{\'e}ter},
  title = {{Linear Algebra Methods in Combinatorics}},
  edition = {Version 2.1},
  publisher = {Department of Computer Science, University of Chicago},
  year = {2020},
  addendum = {\url{https://people.cs.uchicago.edu/~laci/CLASS/HANDOUTS-COMB/BaFrNew.pdf}}
}

@article{ben-eliezer-hod-lovett-2012-low-degree-polys,
  author = {Ben-Eliezer, Ido and Hod, Rani and Lovett, Shachar},
  title = {{Random low-degree polynomials are hard to approximate}},
  journal = {Computational Complexity},
  year = {2012},
  volume = {21},
  number = {1},
  pages = {63--81},
  doi = {10.1007/s00037-011-0020-6},
  addendum = {\url{https://doi.org/10.1007/s00037-011-0020-6}}
}

@article{bernasconi-egidi-hilbert,
  author = {Bernasconi, Anna and Egidi, Lavinia},
  title = {{Hilbert Function and Complexity Lower Bounds for Symmetric Boolean Functions}},
  journal = {Information and Computation},
  year = {1999},
  volume = {153},
  number = {1},
  pages = {1--25},
  doi = {10.1006/inco.1999.2798},
  addendum = {\url{https://doi.org/10.1006/inco.1999.2798}}
}

@article{canfield1980sperner,
  author = {Canfield, E. Rodney},
  title = {{A Sperner property preserved by product}},
  journal = {Linear and Multilinear Algebra},
  year = {1980},
  volume = {9},
  number = {2},
  pages = {151--157},
  doi = {10.1080/03081088008817361},
  addendum = {\url{https://doi.org/10.1080/03081088008817361}}
}

@book{cox2015ideals,
  author = {Cox, David A. and Little, John and O'Shea, Donal},
  title = {{Ideals, Varieties, and Algorithms}},
  edition = {Fourth},
  publisher = {Springer},
  year = {2015},
  doi = {10.1007/978-3-319-16721-3},
  addendum = {\url{https://doi.org/10.1007/978-3-319-16721-3}}
}

@article{debruijn1951set,
  author = {de Bruijn, N. G. and van Ebbenhorst Tengbergen, C. and Kruyswijk, D.},
  title = {{On the set of divisors of a number}},
  journal = {Nieuw Archief voor Wiskunde},
  series = {2},
  year = {1951},
  volume = {23},
  pages = {191--193},
  addendum = {\url{https://research.tue.nl/files/4373475/597494.pdf}}
}

@article{dhand-2014-grid-unimodal,
  author = {Dhand, Vivek},
  title = {{A combinatorial proof of strict unimodality for \(q\)-binomial coefficients}},
  journal = {Discrete Mathematics},
  year = {2014},
  volume = {335},
  pages = {20--24},
  doi = {10.1016/j.disc.2014.07.001},
  addendum = {\url{https://doi.org/10.1016/j.disc.2014.07.001}}
}

@article{felszeghy-hegedus-ronyai-2009-complete-wide,
  author = {Felszeghy, B{\'a}lint and Heged{\H{u}}s, G{\'a}bor and R{\'o}nyai, Lajos},
  title = {{Algebraic Properties of Modulo \(q\) Complete \(\ell\)-Wide Families}},
  journal = {Combinatorics, Probability and Computing},
  year = {2009},
  volume = {18},
  number = {3},
  pages = {309--333},
  doi = {10.1017/S0963548308009619},
  addendum = {\url{https://doi.org/10.1017/S0963548308009619}}
}

@article{hegedus-ronyai-2003-grobner-complete-uniform,
  author = {Heged{\H{u}}s, G{\'a}bor and R{\'o}nyai, Lajos},
  title = {{Gr{\"o}bner bases for complete uniform families}},
  journal = {Journal of Algebraic Combinatorics},
  year = {2003},
  volume = {17},
  number = {2},
  pages = {171--180},
  doi = {10.1023/A:1022934815185},
  addendum = {\url{https://doi.org/10.1023/A:1022934815185}}
}

@article{hegedus-ronyai-2018-linear-sperner,
  author = {Heged{\H{u}}s, G{\'a}bor and R{\'o}nyai, Lajos},
  title = {{A note on linear Sperner families}},
  journal = {Algebra Universalis},
  year = {2018},
  volume = {79},
  eid = {2},
  pages = {1--10},
  doi = {10.1007/s00012-018-0482-3},
  addendum = {\url{https://doi.org/10.1007/s00012-018-0482-3}}
}

@article{heijnen-pellikaan-1998-GHM-Reed-Muller,
  author = {Heijnen, P. and Pellikaan, R.},
  title = {{Generalized Hamming weights of \(q\)-ary Reed--Muller codes}},
  journal = {IEEE Transactions on Information Theory},
  year = {1998},
  volume = {44},
  number = {1},
  pages = {181--196},
  doi = {10.1109/18.651015},
  addendum = {\url{https://doi.org/10.1109/18.651015}}
}

@article{keevash-sudakov-2005-min-rank-inclusion,
  author = {Keevash, Peter and Sudakov, Benny},
  title = {{Set systems with restricted cross-intersections and the minimum rank of inclusion matrices}},
  journal = {SIAM Journal on Discrete Mathematics},
  year = {2005},
  volume = {18},
  number = {4},
  pages = {713--727},
  doi = {10.1137/S0895480103434634},
  addendum = {\url{https://doi.org/10.1137/S0895480103434634}}
}

@article{nie2015hilbert,
  author = {Nie, Zipei and Wang, Anthony Y.},
  title = {{Hilbert functions and the finite degree Zariski closure in finite field combinatorial geometry}},
  journal = {Journal of Combinatorial Theory, Series A},
  year = {2015},
  volume = {134},
  pages = {196--220},
  doi = {10.1016/j.jcta.2015.03.011},
  addendum = {\url{https://doi.org/10.1016/j.jcta.2015.03.011}}
}

@article{proctor1980product,
  author = {Proctor, Robert A. and Saks, Michael E. and Sturtevant, Dean G.},
  title = {{Product partial orders with the Sperner property}},
  journal = {Discrete Mathematics},
  year = {1980},
  volume = {30},
  number = {2},
  pages = {173--180},
  doi = {10.1016/0012-365X(80)90118-1},
  addendum = {\url{https://doi.org/10.1016/0012-365X(80)90118-1}}
}

@article{proctor1982representations,
  author = {Proctor, Robert A.},
  title = {{Representations of \(\mathrm{SL}(2,\mathbb{C})\) on posets and the Sperner property}},
  journal = {SIAM Journal on Algebraic Discrete Methods},
  year = {1982},
  volume = {3},
  number = {2},
  pages = {275--280},
  doi = {10.1137/0603026},
  addendum = {\url{https://doi.org/10.1137/0603026}}
}

@inproceedings{smolensky-low-degree,
  author = {Smolensky, Roman},
  title = {{On representations by low-degree polynomials}},
  booktitle = {Proceedings of the 34th Annual IEEE Symposium on Foundations of Computer Science},
  year = {1993},
  pages = {130--138},
  doi = {10.1109/SFCS.1993.366874},
  addendum = {\url{https://doi.org/10.1109/SFCS.1993.366874}}
}

@article{srinivasan2011symmetric,
  author = {Srinivasan, Murali K.},
  title = {{Symmetric chains, Gelfand--Tsetlin chains, and the Terwilliger algebra of the binary Hamming scheme}},
  journal = {Journal of Algebraic Combinatorics},
  year = {2011},
  volume = {34},
  number = {2},
  pages = {301--322},
  doi = {10.1007/s10801-010-0272-2},
  addendum = {\url{https://doi.org/10.1007/s10801-010-0272-2}}
}

@online{venkitesh-2021-zariski-positive-char,
  author = {Srinivasan, Srikanth and Venkitesh, S.},
  title = {{On Vanishing Properties of Polynomials on Symmetric Sets of the Boolean Cube, in Positive Characteristic}},
  year = {2021},
  eprinttype = {arXiv},
  eprint = {2111.05445},
  addendum = {\url{https://arxiv.org/abs/2111.05445}}
}

@article{venkitesh-2022-covering-symmetric-sets,
  author = {Venkitesh, S.},
  title = {{Covering Symmetric Sets of the Boolean Cube by Affine Hyperplanes}},
  journal = {The Electronic Journal of Combinatorics},
  year = {2022},
  volume = {29},
  number = {2},
  eid = {P2.22},
  doi = {10.37236/10600},
  addendum = {\url{https://doi.org/10.37236/10600}}
}

@article{wei-1991-GHM,
  author = {Wei, Victor K.},
  title = {{Generalized Hamming Weights for Linear Codes}},
  journal = {IEEE Transactions on Information Theory},
  year = {1991},
  volume = {37},
  number = {5},
  pages = {1412--1418},
  doi = {10.1109/18.133259},
  addendum = {\url{https://doi.org/10.1109/18.133259}}
}

@article{wilson-1990-diagonal-incidence,
  author = {Wilson, Richard M.},
  title = {{A diagonal form for the incidence matrices of \(t\)-subsets vs. \(k\)-subsets}},
  journal = {European Journal of Combinatorics},
  year = {1990},
  volume = {11},
  number = {6},
  pages = {609--615},
  doi = {10.1016/S0195-6698(13)80046-7},
  addendum = {\url{https://doi.org/10.1016/S0195-6698(13)80046-7}}
}

@inproceedings{GGHNY-2024-hilbert,
  author = {Golovnev, Alexander and Guo, Zeyu and Hatami, Pooya and Nagargoje, Satyajeet and Yan, Chao},
  title = {{Hilbert Functions and Low-Degree Randomness Extractors}},
  booktitle = {Approximation, Randomization, and Combinatorial Optimization. Algorithms and Techniques (APPROX/RANDOM 2024)},
  series = {Leibniz International Proceedings in Informatics (LIPIcs)},
  volume = {317},
  pages = {41:1--41:24},
  publisher = {Schloss Dagstuhl -- Leibniz-Zentrum f{\"u}r Informatik},
  year = {2024},
  doi = {10.4230/LIPIcs.APPROX/RANDOM.2024.41},
  addendum = {\url{https://doi.org/10.4230/LIPIcs.APPROX/RANDOM.2024.41}}
}
\cleardoublepage
\restoregeometry

\end{document}